\documentclass[a4paper,11pt]{amsart}
\usepackage{amssymb}
\usepackage[T1]{fontenc}
\usepackage[ansinew]{inputenc}

\usepackage{hyperref}

\usepackage{amsrefs}

\newtheorem{theorem}{Theorem}[section]
\newtheorem{proposition}{Proposition}[section]
\newtheorem{corollary}{Corollary}[section]
\newtheorem{lemma}{Lemma}[section]
\newtheorem{remark}{Remark}[section]
\newtheorem{example}{Example}[section]
\newtheorem{definition}{Definition}[section]

\newcommand{\0}{{\boldsymbol{0}}}
\newcommand{\1}{{\boldsymbol{1}}}

\newcommand{\N}{{\mathbb N}}
\newcommand{\Q}{{\mathbb Q}}
\newcommand{\R}{{\mathbb R}}
\newcommand{\Z}{{\mathbb Z}}
\newcommand{\K}{{\mathbb K}}

\newcommand{\ad}{\operatorname{ad}}
\newcommand{\Ber}{\operatorname{Ber}}
\newcommand{\diverg}{\operatorname{div}}
\newcommand{\End}{\operatorname{End}}
\newcommand{\GL}{\operatorname{GL}}
\newcommand{\Hom}{\operatorname{Hom}}
\newcommand{\id}{\operatorname{id}}
\newcommand{\im}{\operatorname{im}}
\newcommand{\sh}{\operatorname{sh}}
\newcommand{\str}{\operatorname{str}}
\renewcommand{\th}{\operatorname{th}}
\newcommand{\triv}{\operatorname{triv}}

\newcommand{\CL}{{\mathcal L}}
\newcommand{\CV}{{\mathcal V}}
\newcommand{\CT}{{\mathcal T}}
\newcommand{\CW}{{\mathcal W}}

\renewcommand{\a}{\mathfrak{a}}
\newcommand{\g}{\mathfrak{g}}
\newcommand{\h}{\mathfrak{h}}
\newcommand{\osp}{{\mathfrak o\mathfrak s\mathfrak p}}
\newcommand{\p}{{\mathfrak p}}
\newcommand{\q}{{\mathfrak q}}

\begin{document}

\title{Symmetric pairs and Gorelik elements}
\dedicatory{To Ernest Borisovich  Vinberg on the occasion of his
seventieth birthday}
\author{Michel DUFLO}
\address{University Denis Diderot-Paris~7,  Institut de Math\'ematiques de Jussieu,
 C.P.~7012\\ 2~place Jussieu,   F-75251 Paris~cedex~05}
\email{duflo@math.jussieu.fr}
\author{Emanuela PETRACCI}
\address{Department of Mathematics, University of Toronto, 40 St. George St.
Toronto, Ontario, Canada M5S 2E4}

\email{petracci@math.toronto.edu}
\thanks{The  research of the second author was done while she was working
at  the Mathematics department {\it Section de math\'ematiques} of
the University of Geneva and received some financial support from
the Swiss National Science Foundation. During  part of the final
redaction she  was   supported first by   a  post-doctoral grant
of the European  network LIEGRITS at the I.E.C.N. Institute in
Nancy, and then by an invitation at the I.H.E.S Institute in
Bures-sur-Yvettes. The second author thanks   all these
institutions. A special thanks to Prof.~Caroline Gruson, LIEGRITS
responsible  in Nancy.}

\begin{abstract}
The second author gave a formula for   the elements of the
enveloping algebra of a  Lie superalgebra  defined by  Gorelik
under an appropriate unimodularity assumption. We show that this
formula is a particular case of a formula for the Jacobian of the
exponential map of a symmetric superspace.
\end{abstract}



\begin{center}
\begin{tabular}{l}
\scriptsize{To appear in the special volume of Journal of Algebra 
dedicated to E.B. Vinberg}
\\
\hline
\end{tabular}
\end{center}

\maketitle
{\scriptsize \tableofcontents}
\section*{Introduction.}

In her PhD thesis \cite{Petracci-thesis}, the second author gave a
formula for the  elements, in the enveloping algebra of a Lie
superalgebra,  defined by  Gorelik  \cite{Gorelik2001} under an
appropriate unimodularity assumption. We recall Gorelik's
definition. Let $\g=\g_\0\oplus \g_\1$ be a finite dimensional
superalgebra over a field of characteristic $0$, and $\sigma$ be
the automorphism of $\g$ such that $\sigma (a)=a$ if $a\in \g_\0$
and $\sigma (a)=-a$ if $a\in \g_\1$. Arnaudon, Bauer and Frappat
\cite{ABF1997} define the \emph{twisted adjoint action} $\ad'(a)$
of $a\in \g$ in the enveloping algebra $U(\g)$  by
\begin{equation}\label{eq:twisted}
\ad'( a)(u)=au- \pm u \sigma(a)  \quad\text {for }\quad  u \in
U(\g).
\end{equation}
Let $S(\g)$ be the symmetric algebra and $\beta :S(\g)\to U(\g)$
be the symmetrization map.  The finite dimensional\footnote{Since
$\g_\1$ is purely odd, $S(\g_\1)$ is isomorphic to the exterior
algebra of $\g_\1$ (see remark \ref{rem:exterior}).} subspace
$\beta(S(\g_\1))$ of $U(\g)$ is stable under the twisted adjoint
action (see~\cite{ABF1997}). Gorelik \cite{Gorelik2001} proves
that \emph{if the following unimodularity condition holds},
\begin{equation}\label{eq:unim}
\str_{\g_\1}(\ad a)=0 \quad\text {for }\quad  a \in \g_\0,
\end{equation}
there exists a non zero invariant of the twisted adjoint action $T
\in \beta(S(\g_\1))$, unique up to a multiplicative constant. We
 call $T$ a \emph{Gorelik element} of $U(\g)$. The Gorelik
elements occurred already for $\osp(1,2)$ in
\citelist{\cite{Pinczon1990} \cite{Lesniewski1995}},  for
$\osp(1,2l)$ in \cite{ABF1997} under the name \emph{Casimir
ghosts}, and in \cite{Musson1997}. The Gorelik elements play an
important role in the representation theory of $\g$, see in
particular \citelist{ \cite{Pinczon1990} \cite{Musson1997}
\cite{GL1999}  \cite{Gorelik2001} \cite{Gorelik2002}}.

\medskip
\noindent {\bf Formula for Gorelik's elements.} Let us describe
the formula given in  \cite{Petracci-thesis} for the Gorelik
elements. Let $q=\dim\g_1$. Choose a basis $d$ of the
$1$-dimensional space $S^q(\g_\1)$.  Consider the algebra
$S(\g_\1)^*$ dual of the coalgebra $S(\g_\1)$. By duality,
$S(\g_\1)$ is a module over $S(\g_\1)^*$; this action is the usual
interior product and the map $f \to fd$ is an isomorphism from
$S(\g_\1)^*$ to $S(\g_\1)$. In particular, to  describe our
formula means to describe a special element $J_2\in S(\g_\1)^*$.
To do this, we use the \emph{generic point  of $\g_1$}. It is the
element $y=\sum e_i x^i $ of the Lie $S(\g_\1)^*$-superalgebra
$\g\otimes S(\g_\1)^*$, where $(e_i)_{1\leq i \leq q}$ is a basis
of $\g_1$ and $(x^i)_{1\leq i \leq q}$ is the corresponding dual
basis of $\g_1^*$. For $k\in \N$, we have
\begin{equation}\label{eq:adyk}
  \ad^{k} y
= \sum_{1\leq i_1<i_2<\dots <i_k\leq q\ }
  \sum_{s\in S_k}
  \pm
  \ad e_{i_{s(1)}} \cdots \ad e_{i_{s(k)}}\  x^{i_1}\cdots x^{i_k}.
\end{equation}
In particular, $\ad^k y=0$ if $k>q$. When   $k$ is even the space
$\g_\1$ is invariant by $\ad e_{i_{s(1)}} \cdots \ad
e_{i_{s(k)}}$. We consider the finite sum
\begin{equation}\label{eq:bersh}
 \frac{\sh(\ad\frac{ y}{2})}{\ad\frac{ y}{2}}
 = 1+ \frac{\ad^2 y}{24} +\frac{\ad^4 y}{1920}+\cdots
\end{equation}
It is an even invertible endomorphism of $\g\otimes S(\g_\1)^*$
which stabilizes $\g_\1\otimes S(\g_\1)^*$. We denote by   $J_2
\in S(\g_\1)^*$ the Berezinian (in this case the inverse of the
determinant) of this endomorphism of $\g_\1\otimes S(\g_\1)^*$. We
summarize this definition by the notation
\begin{equation}\label{eq:J}
J_2= \Ber_{\g_1}\left(\frac{\sh(\ad \frac{ y}{2})}{\ad\frac{
y}{2}}\right).
\end{equation}
Similarly, if $k$ is even, we define
\begin{equation}\label{eq:adyk2}
\str_{\g_1}(\ad^{k} y)
=\sum_{1\leq i_1<i_2<\dots <i_k\leq q\ }\sum_{s\in S_k}
\pm \str_{\g_1}(\ad e_{i_{s1}} \cdots \ad e_{i_{sk}})\
x^{i_1}\cdots x^{i_k} \in S(\g_1)^*.
\end{equation}
The relation   $\Ber( \exp(\cdot))= \exp(\str(\cdot))$ allows to
write explicitly the Taylor expansion of $J_2$ in  homogeneous
components:
\begin{equation}\label{eq:Ja}
J_2=1 +   \frac{1}{24} \str_{\g_1}(\ad^{2} y)
-\frac{\str_{\g_1}(\ad^{4} y)}{2880}+ \frac{ \str^2_{\g_1}(
\ad^{2} y )}{1152} + \cdots
\end{equation}
For instance, for $q=2$ we obtain
\begin{equation}\label{eq:Jc}
J_2= 1 +   \frac{1}{24} \str_{\g_1}(-\ad e_1 \ad e_2 +\ad e_2  \ad
e_1)\ x^1 x^2.
\end{equation}

We have ---see \cite{Petracci-thesis} and Corollary
\ref{cor:Gorelikfinal} below:
\begin{theorem}\label{theo:J}
Under the unimodularity condition (\ref{eq:unim}), $\beta(J_2 d)$
is a Gorelik element of $U(\g)$.
\end{theorem}

For instance, for $q=2$, if  we   choose $d=e_1e_2$ we obtain the
Gorelik element
\begin{equation}\label{eq:Tb}
T= \frac{1}{2}(j(e_1)j(e_2)-j(e_2)j(e_1)) + \frac{1}{24}
\str_{\g_1}(\ad e_1  \ad e_2 -\ad e_2  \ad e_1 ) ,
\end{equation}
where $j: \g \to U(\g)$ is the canonical injection.

\medskip
\noindent {\bf Occurrence of formal symmetric spaces.} The
occurrence of the  operator $\frac{\sh(\ad\frac{ y}{2})}{\ad
\frac{ y}{2}}$ suggests a relation with  the differential of the
exponential map. Indeed, \emph{we prove  in this paper that
Theorem \ref{theo:J} is  a particular case of a formula for the
Jacobian of the exponential map for symmetric spaces ---in the
sense of supergeometry}. A natural setting to present this result
is to consider a Lie superalgebra $\g$ over a commutative
$\Q$-superalgebra $\K$, and a subalgebra $\h \subset \g$ such that
the quotient $\g/\h$ is free of finite rank
 over $\K$. The formal Lie supergroups $G$, $H$
and the formal homogeneous space $G/H$ can be defined  by mean of
suitable commutative $\K$-superalgebras $F(G)$, $F(H)$, $F(G/H)$
of formal functions. Formal supergroups were first considered by
Berezin and Kac \cite{Berezin1970}. There exists and we choose a
formal diffeomorphism  $\theta: \g/\h \to G/H$, or equivalently an
algebra isomorphism $\theta^*:   F(G/H) \to F(\g/\h):=S(\g/\h)^*$
--- i.e.  into
  the algebra
 of formal ``functions'' on $\g/\h$, which
is an algebra of formal power series in a finite number of even
and odd variables.   In this situation, the \emph{dualizing
$F(G/H)$-module}  $K(F(G/H))$ is defined. It is free of rank $1$:
in the purely even case (i.e. $\K=\K_\0$ and  $\g=\g_\0$) it is
the module of top degree differential forms on $G/H$, and in
general it is the Berezinian of the module of $1$-differential
forms (see \cite{MPV1990} and subsection \ref{sec:canonical}
below). There is a natural action of $\g$ in $K(F(G/H))$ by Lie
derivatives, and a natural $\K$-linear isomorphism
$\theta^*:K(F(G/H))\to K(F(\g/\h))$.  The unimodularity assumption
reads:
\begin{equation}\label{eq:unimb}
\str_{\g/\h}(\ad a)=0 \quad\text {for }\quad  a \in \h.
\end{equation}
We show in proposition \ref{prop:invariant} that condition
(\ref{eq:unimb}) is equivalent to the existence of a
$\g$-invariant basis $D_{G/H}\in K(F(G/H))$. In this case, there
is a translation invariant element $d_{\g/\h}\in K(F(\g/\h))$ and
a function $J\in F(\g/\h)$, the \emph{Jacobian} of $\theta$, such
that
\begin{equation}\label{eq:jacobian}
\theta^*(D_{G/H})=J d_{\g/\h}  \mbox{ and } J(0)=1.
\end{equation}
We compute this Jacobian in two cases.

\medskip
First,  we suppose that  $\h=0$ and that $\g$ is a free
$\K$-module of finite rank. In theorem \ref{theo:jac} we give  a
formula for the Jacobian of the exponential map $\g \to G$. The
result is not surprising, but we think that it is a good
introduction   to the second case.

\medskip
In the second case, we assume  that $\h$ is the fixed point set of
an involution $\sigma$ of $\g$. We write $\g=\h\oplus \q$ the
decomposition of $\g$ in $1$ and $-1$ eigenspaces of $\sigma$. We
say that $(\g,\sigma)$ or $\g=\h\oplus \q$  is a symmetric pair,
and that $G/H$ is symmetric. The exponential map is a formal
diffeomorphism $\q \to G/H$. We prove (see theorem \ref{thm:main}
below):
\begin{theorem}\label{theo:Ja}
Suppose that $\g=\h\oplus \q$ is a symmetric pair, that $\q$ is
free of finite rank, and that  condition (\ref{eq:unimb}) holds.
Then the Jacobian of the exponential map $\q \to G/H$ is the
element $J \in F(\q)$ given by the formula
\begin{equation}\label{eq:jacobian4}
J= \Ber_{\q}\left(\frac{\sh(\ad  y)}{\ad y}\right),
\end{equation}
where   $y$ is the generic point of $\q$.
\end{theorem}

\medskip
We illustrate Theorem \ref{theo:Ja} in two particular extreme
cases.
\begin{example}
Suppose that $\K=\R$ and that $\g=\g_\0$ is a finite dimensional
Lie algebra. Let $\mathbf G$  and  $\mathbf H$ be connected Lie
groups with Lie algebras $\g$ and $\h$, and such that $\mathbf H$
is a closed subgroup of $\mathbf G$. Condition (\ref{eq:unimb}) is
equivalent to the existence of a non zero left invariant
differential form of maximal degree $D_{\mathbf G/\mathbf H}$ on
$\mathbf G/\mathbf H$. It is well known (see for instance
\cite{Rouviere1986}) that its inverse image by the exponential map
$\exp: \q \to \mathbf G/\mathbf H$ at a point $a\in \q$ is equal
to $\det_{\q}(\frac{\sh(\ad a)}{\ad a}) d_{\q}$, where $d_{\q}$ is
a form of maximal degree  on $\q$ invariant by translations. Since
in this case $\Ber_\q=\det_\q$, we recover Theorem \ref{theo:Ja}
by taking the Taylor series of $\det_{\q}(\frac{\sh(\ad a)}{\ad
a})$ at $a=0$.
\end{example}
\begin{example}
Suppose that $\K=\K_\0$ is a field of characteristic $0$, and
consider a finite dimensional Lie superalgebra $\g=\g_\0\oplus
\g_\1$. This is a symmetric pair, and   condition (\ref{eq:unimb})
is the same than   Gorelik's condition (\ref{eq:unim}). Assuming
this condition,  we prove in  theorem \ref{thm:correspondance} a
result on induced representations of a Lie superalgebra which says
in particular that \emph{the image $\widetilde \beta(J d)$ of $
\beta(J d)$ in $U(\g)/U(\g)\g_\0$ is a non zero $\g$-invariant
element}. We explain in subsection  \ref{sec:gorelik} how Theorem
\ref{theo:J} follows  by considering an isomorphism of
$\g$-modules from $U(\g)/U(\g)\g_\0$ to the subspace
$\beta(S(\g_\1)) $ of $U(\g)$ endowed with the twisted adjoint
action.
\end{example}

\medskip
In fact, it is better not to assume  condition (\ref{eq:unimb}):
in general, we twist $K(F(G/H))$ by a suitable line bundle on
$G/H$ to get the existence of an invariant element. More
precisely, we consider the $\g$-module $F(G/H, \Ber{\g/\h})$
coinduced by the rank $1$ $\h$-module provided by the character
$\str_{\g/\h}$ (for the details of these notations, see section
\ref{sec:inducedformal}). It is a free $F(G/H)$-module of rank
$1$, as well as the \emph{twisted dualizing module} $
K(F(G/H),\Ber{\g/\h}):=K(F(G/H))\otimes_{F(G/H)}F(G/H,
\Ber{\g/\h})$. We prove in proposition \ref{prop:canonical} that
there is a canonical $\g$-invariant basis $ \mathbf D_{G/H}\in
K(F(G/H),\Ber{\g/\h})$. In Theorem \ref{thm:main},  we  prove
formula (\ref{eq:jacobian4}) also in the twisted case.

\medskip
\noindent {\bf Formulas related to the symmetrization map.} Our
proofs involve many useful formulas about the symmetrization
$\beta$, in particular how to work with   induced and coinduced
representations in exponential coordinates. We think that they are
of interest, independently of the present application to
Jacobians. In the Lie group case $\mathbf G/\mathbf H$  these
formulas occur more or less explicitly in Rouvi{\`e}re
\citelist{\cite{Rouviere1986}\cite{Rouviere1994}}. In the case of
the symmetric pair $\g=\g_\0\oplus \g_\1$ associated to a finite
dimensional Lie superalgebra defined over a field of
characteristic zero, they are given in Koszul \cite{Koszul1982}.
However, our proofs are different. We use the methods of
\citelist{\cite{Petracci-thesis}  \cite{Petracci2003}}: we seek
\emph{universal  formulas}, that is formulas which can be
expressed in the same manner for all symmetric pairs $ \g=\h\oplus
\q$. They depend on a finite number of formal power series in
$\K_\0[[t]]$ which must satisfy some \emph{functional equations}.
We are lucky enough: these functional equations have solutions
which are essentially unique, and which are given by explicit
formulas.

\medskip
\noindent {\bf Hypotheses used in the text.}  We assume from the
beginning  that $ \K$ contains   $\Q$, because the main actors of
this paper are the symmetrization $\beta : S(\g) \to U(\g)$ and
the related exponential map $\g \to G$. But nothing more. Besides
the obvious fact that it is better to avoid unnecessary
hypotheses, adding more assumptions on $\K$ (like $\K$ being a
field, or the weaker assumption $\K=\K_\0$) would not be specially
helpful, because, as we have seen in formula (\ref{eq:J}), it is
convenient to use Lie superalgebras over rings like exterior
algebras.

\medskip
For simplicity, we state and prove our results involving
Berezinians and friends when $\g/\h$ is free of finite rank, even
if they can probably easily be extended to the case when $\g/\h$
is projective of finite rank, or to a locally free sheaf version
of it. For the other results, we suppose that $\h$ is a direct
factor of $\g$, but we do not require anything more on the
$\K$-module structure of $\g$ or $\h$.

\medskip
Working in this generality forces us to write more natural proofs.
The drawback is that we have to use  some results usually stated
under the assumption that $\K$ is a field. The main one is the
Poincar{\'e}-Birkhoff-Witt theorem which states that the
symmetrization $\beta : S(\g) \to U(\g)$  is an isomorphism of
coalgebras. This one is proved  in \cite{Cohn1963} --see also
\cite{Petracci2003}. For some others, like Blattner's
imprimitivity theorem \cite{Blattner1969}, we provide the proofs.

\section{Berezinians and dualizing modules.}
\label{SEC:1}
\subsection{Notations for modules.}\label{sub:0}
In this text, \emph{all $\Z$-modules are supposed to be
$\Z/2\Z$-graded}. Let $M=M_\0 \oplus M_\1$ be a $\Z$-module. The
elements of $M_\0$ are called even, and the elements of $M_\1$
odd. We denote by $p(m)\in \Z/2\Z=\{\0,\1\}$ the parity of an
element $m$, in the following sense: whenever the expression
$p(m)$ occurs in a formula, it is implicitly assumed that $m \in
M_{p(m)}$.

Depending on the context, we use  the notations $\id_M$, $\id$ or
$1$ for the identity map of $M$.  We write $\Pi $ for the change
of parity functor. Thus $(\Pi M)_\0=M_\1$, $(\Pi M)_\1=M_\0$. We
use the notations  $\pi_M$ or $\pi$ for the identity map of $M$
when it is considered as a map from $M \to \Pi M$.

We fix a   commutative  $\Q$-superalgebra $\K=\K_\0\oplus\K_\1 $,
associative with unit $1\in \K_\0$  (supercommutative in the
terminology of \cite{Manin1997}*{ch. 3}).

If $M=M_\0 \oplus M_\1$ is a $\K$-module, we write $km=\pm mk$ for
$k \in \K$ and $m\in M$, \emph{where here and after, the symbols
are implicitly assumed to be homogeneous, and the sign $\pm$ is
obtained by the standard rule of signs} (see \cite{Manin1997}*{ch.
3.1}). For instance in this case, $\pm=(-1)^{p(k)p(m)}$.

If $M$ and $N$ are  $\K$-modules, if $\K$ is clear from the
context we write $M\otimes N$ for $M\otimes_\K N$, $\Hom(M,N)$ for
$\Hom_\K(M,N)$, $\End(M)$ for $\Hom(M,M)$, $M^*$ for $\Hom(M,\K)$.
We recall that the $\K$-structure on $\Pi M$ is defined by
$\pi_M(m)k=\pi_M(mk)$, i.e. by the requirement $\pi_M \in
\Hom(M,\Pi M)_\1$.

Let   $M$ be a  free $\K$-module of finite rank. It is isomorphic
to $\K^p\oplus \Pi\K^q$ for some uniquely determined $p$ and $q$
in $\N$. We say that  $p+ q \in \N$ is the rank of $M$. If we need
to keep track of parities, we say that the rank is $(p,q)\in
\N^2$. For instance, a module of rank $1$ is isomorphic either to
$\K$ (in this case it is of rank $(1,0)$) or to $\Pi \K$ (in this
case it is of rank $(0,1)$).

Suppose $(e_i)_{i\in I}$ is a basis of a free module $M$ of finite
rank. {\it In this text we always suppose that basis are
 homogeneous}. We
denote by $(x^i)_{i\in I}$ the dual basis of $M^*$:
$x^i(e_j)=\delta^i_j$. If $N$ is a $\K$-module, we identify
$\Hom(M,N)$ and $N\otimes M^*$ in the usual manner. This
identification sends $\id_M \in \End(M)_\0$ to $\sum e_i\otimes
x^i$, and $\pi_M \in \Hom(M,\Pi M)_\1$ to $\sum \pi_M e_i\otimes
x^i$.

\subsection{Symmetric algebras as coalgebras.}
\label{sub:1} In this subsection, we state the properties of
formal functions on a $\K$-module which we shall need in the
sequel. For a more detailed exposition, we refer to
\citelist{\cite{Radford1986} \cite{Petracci2003}}. Let $M$ be a
$\K$-module. We denote by $S(M)=\oplus_{n\geq 0} S^n(M)$ the
symmetric algebra of $M$. It is an Hopf superalgebra.

\begin{remark}\label{rem:exterior}
Everything in this article is  considered in the super setting.
Thus, $S(M)$ is the commutative superalgebra freely generated by
$M.$ If $M$ is a free $\K$-module with basis $(e_i)_{1\leq i \leq
n}$, $S(M)$ is equal to the commutative superalgebra
$\K[e_1,\dots,e_n]$ of finite sums with coefficients in $\K$ of
monomials $e_1^{p_1}\dots e_n^{p_n}$, with $p_i\leq 1$ if $e_i$ is
odd.

In particular, if $(e_1, ...,e_q)$ is a basis of $M$ consisting of
$q$ odd elements,  $S(M)=\K[e_1,\dots,e_q]$ is the exterior
algebra on the generators $e_i$. It is a free $\K$-module of rank
$2^q$.
\end{remark}
We denote  by $\Delta_{S(M)}\in \Hom(S(M),S(M)\otimes S(M))$ the
coproduct, and by $S_{S(M)}: S(M)\to S(M)$ the  antipode. Let
$w\in S(M)$, we denote its coproduct by
$$
\Delta_{ S(M)}(w)=\sum w_i\otimes w'_i .$$

We denote by $F(M)$ the algebra $S(M)^*$ dual to the coalgebra
$S(M)$, and we consider it as the algebra of formal $\K$-valued
functions defined in a formal neighborhood of $0\in M$. We denote
by $\delta: F(M)\to \K$ the map defined by $\delta(f)=f(1)$. It is
a morphism of $\K$-superalgebras, and we consider it as the
``evaluation at $0\in M$''. Let $S^+(M)=\ker \delta$. Thus
$S^+(M)=\oplus_{n>0} S^n(M)$.

Let $N$ be a $\K$-module. We denote by $F(M,N)= \Hom(S(M),N)$ the
space of formal $N$-valued functions on $M$.  We identify  $N$
with the submodule of $F(M,N)$ consisting of the ``constant
functions'', that is the functions which are zero over $S^+(M)$.
We still denote by $\delta: F(M,N)\to N$ the map defined by
$\delta(f)=f(1)$.

\medskip
The space   $F(M,N)$ is an $F(M)$-module. More generally, if $\mu:
N_1\otimes N_2 \to N$ is a map of $\K$-modules,  there is a
corresponding map $\tilde\mu: F(M,N_1)\otimes_{F(M)} F(M,N_2)\to
F(M,N)$, defined by the formula
\begin{equation}\label{eq:0}
\tilde \mu(f_1\otimes f_2)=\mu\circ (f_1\otimes f_2)
\circ\Delta_{S(M)} , \ \ {\rm for}\  f_1\in F(M,N_1), f_2\in
F(M,N_2) .
\end{equation}
Depending on the setting, there are various notations for $\tilde
\mu(f_1\otimes f_2)$. For instance, if $N$ is a Lie
$\K$-superalgebra with product $[n_1,n_2]=\mu(n_1\otimes n_2)$, we
usually employ  the notation $[f_1,f_2]=\tilde \mu(f_1\otimes
f_2)$, which gives to $F(M,N)$ a Lie $F(M)$-superalgebra
structure. Similarly, if $N$ is an associative algebra with
product $n_1n_2= \mu(n_1\otimes n_2)$, we usually employ  the
notation $f_1f_2=\tilde \mu(f_1\otimes f_2)$.

However, in the special case $N=S(M)$ we also use the notation
$f_1*f_2=\tilde \mu(f_1\otimes f_2)$. Indeed, it is the
traditional notation (see e.g. \cite{Radford1986}), and it helps
to avoid   confusions with the   composition of maps  in $
F(M,S(M))=\End(S(M))$, which is a different associative product.
The constant function $1 \in F(M,S(M))$ is the unit for the
product $*$. The identity $\id_{S(M)} \in \End_\K(S(M))$ and the
antipode $S_{S(M)} \in \End_\K(S(M))$  are mutual inverses for the
product $*$.

\medskip
Let $\phi \in \End_\K(S(M))$. We  use the notation $c_\phi=\phi
*\id_{S(M)}$ and we denote by $\zeta_\phi$ its transpose acting
(as a derivation) in $F(M,N)$. For $w\in S(M)$   and $g\in
F(M,N)$,  we have:
\begin{equation}\label{eq:6}
c_\phi(w)=\sum \phi(w_i)w'_i \mbox{    and   } \zeta_\phi(g)=\pm
g\circ c_\phi.
\end{equation}
Note that any $\psi\in \End_\K(S(M))$ is of the form $\psi=c_\phi$
with $\phi=S_{S(M)}*\psi$. We consider particular cases.
\begin{enumerate}
\item[i)] Let $f\in F(M)$. Then $\zeta_f\in \End_\K(F(M,N))$  is
the left multiplication by $f$ which defines the $F(M)$-module
structure of $F(M,N)$. Thus we  just write $fg=\zeta_f(g)$ for
$g\in F(M,N)$. Accordingly, for $w\in S(M)$, we also write $fw$
for $c_f(w)$. {\it This provides $S(M)$ with a structure of
$F(M)$-module}.

Note that if $f\in M^*$, then $c_f$ is a derivation of $S(M)$, it
is the
 derivative in the direction $f$. So, for a general $f\in F(M)$,
$c_f$ may be thought of as an infinite order differential operator
with constant coefficients on $S(M)$. However, in this paper, we
prefer to  consider  $c_f(w)$ as  ``the product of the
distribution $w$ by the function $f$''.
\item[ii)] Let $w \in S(M)$. Then $c_w$ is the left multiplication
by $w$. Thus, for $f\in F(M,N)$ and $w'\in S(M)$ we have:
\begin{equation}\label{eq:5}
\zeta_w(f)(w')=\pm f(ww').
\end{equation}
We consider $\zeta_w$ as a constant coefficient differential
operator on $F(M,N)$. In particular, for $m\in M$,  $\zeta_m$ is
the derivative in the direction $m$.
\end{enumerate}
\medskip

The space $F(M,M)\subset \End(S(M))$ is called the space of {\it
formal vector fields} on $M$.   Let $\alpha\in F(M,M)$. In this
case $c_\alpha$ is a coderivation of $S(M)$ and $\zeta_\alpha$   a
derivation of $F(M)$.  For  $\alpha$ and $\beta$ in $F(M,M)$, we
have (see \citelist{\cite{Radford1986} \cite{Petracci2003}}):
\begin{equation}\label{eq:11}
  [c_\alpha , c_\beta ]
=-c_\gamma \mbox{ and } [ \zeta_\alpha , \zeta_\beta
]=\zeta_\gamma \mbox{,   where } \gamma=\zeta_\alpha(\beta)-\pm
\zeta_\beta(\alpha).
\end{equation}

\medskip
We denote by $x_M$ the identity map of $M$, considered as the
element of $F(M,M)$ which is $0$ on $S^n(M)$ for $n\neq 1$. We
call it the \emph{generic point of $M$}.

\subsection{Functional calculus.}
Let $M$ and $N$ be   $\K$-modules. We say that a   family
$(f_i)_{i\in I}$ of elements of $F(M,N)$ is locally finite if, for
all $w\in S(M)$, the number of elements $i\in I$ such that
$f_i(w)\neq 0$ is finite. In this case, we denote by
\begin{equation}\label{eq:sum}
\sum_{i\in I} f_i\in F(M,N)
\end{equation}
the element such that
\begin{equation}\label{eq:suma}(\sum_{i\in I}
f_i)(w)=\sum_{i\in I} f_i(w) \mbox{ for  } w\in S(M).
\end{equation}

\medskip
Let $A$ be an associative $\K$-superalgebra with unit $1\in A_\0$,
and
 $M$ be a $\K$-module. Let us consider an $ X  \in F(M,A)_\0$  such that
$\delta( X  )=0$. For  $w\in S^k(M)$, we have $ X ^{n}(w)=0$ for
$n>k$. Thus, for  $f =\sum_{n\geq 0} f_n t^n \in \K_\0[[t]]$ we
can define $f(X)\in F(M,A)_\0$ by the locally finite sum
\begin{equation}\label{eq:function}
f( X  )=\sum_{n\geq 0} f_n  X  ^n .
\end{equation}
The rules of functional calculus apply: $f \to f( X  )$ is a
$\K_\0$-algebra homomorphism, and $\delta(f( X ))=f_0$. Moreover,
if $g \in \K_\0[[t]]$ verifies $g_0=0$, then $f(g( X  ))=f\circ g
(X)$ where  $f\circ g \in \K_\0[[t]]$ is defined by  $f\circ g
(t)=f(g(t))$. In particular, $(1+ X )^{-1}:=\sum_{n\geq 0}(-1)^n X
^n$, $\log(1+ X  ):= \sum_{n> 0} (-1)^{n+1} \frac{ X  ^n}{n}$ and
$\exp( X ):=\sum_{n\geq 0}\frac{ X  ^n}{n!}$ are defined. We shall
need the following standard supplement to the functional calculus
in subsection \ref{sec:generic}.
\begin{lemma}\label{lem:dercomp}
Let $\alpha \in F(M,M)$, $f(t)\in \K_\0[[t]]$, and $X \in F(M)_\0$
be such that $\delta(X)=0$. Then
$\zeta_\alpha(f(X))=f'(X)\zeta_\alpha(X)$.
\end{lemma}
\begin{proof}
We leave to the reader to check that it reduces to the case
$f=t^n$ and to the fact that $\zeta_\alpha$ is a derivation.
\end{proof}

\begin{lemma}\label{lem:basis}
Suppose that $N$ is a free $\K$-module of finite rank $(p,q)$.
Then $F(M,N)$ is a free $F(M)$-module of finite rank $(p,q)$. Let
$(\phi_i)_{i\in I}$ be a family of elements of $F(M,N)$. Then it
is a basis of $F(M,N)$ if and only if the family
$(\delta(\phi_i))_{i\in I}$ is a basis of $N$.
\end{lemma}
\begin{proof}
The  inclusion $N \subset F(M,N)$ induces an isomorphism
$F(M)\otimes N \simeq F(M,N)$. Let $v_i:=\delta(\phi_i)$.
Obviously, if $(\phi_i)_{i\in I}$ is a basis of $F(M,N)$ then
$(v_i)_{i\in I}$ is a basis of $N$. Conversely, suppose that
$(v_i)_{i\in I}$ is a basis of $N$. By the inclusion $N \subset
F(M,N)$ it is a basis of $F(M,N)$. Write $\phi_j= \sum v_i g_j^i$
with $g_j^i \in F(M)$. The square matrix $g:=(g_j^i)$ is even. The
matrix $\delta(g)$ is the identity $I$. By functional calculus $g$
is invertible with inverse $(I+(g-I))^{-1}=I + (I-g)+(I-g)^2+
\cdots$. Thus $(\phi_i)_{i\in I}$ is a basis of $F(M,N)$.
\end{proof}

\subsection{Berezinians and supertraces.} \label{subsec:ber}
In this subsection, we consider a  $ \K$-module $N$ which is  free
of finite rank.

Inside the commutative superalgebra $F(N,S(\Pi N))$ we consider
the subalgebra $S(\Pi N)\otimes S(N^*)=S(\Pi N \oplus N^*)$, and
we still denote by $\pi$ the element of $S(\Pi N \oplus N^*)$
which corresponds to $\pi \in \Hom(N,\Pi N)$. We have $\pi ^2=0$.
Following Manin \cite{Manin1997}*{ch. 3.4}, we define the
\emph{Berezinian of the module} $N$
\begin{equation}\label{eq:berezinianmodule}
\Ber_\K (N):=\ker \pi/ \im \pi
\end{equation}
as the homology of the multiplication by $\pi$ in $S(\Pi N \oplus
N^*)$. We write $\Ber N$ or $\Ber(N)$ for $\Ber_\K (N)$ if $\K$ is
clear from the context.

Let $(e_i)_{1\leq i \leq n}$ be a  totally ordered basis of $N$,
and $(x^i)_{1\leq i \leq n}$ be the dual basis of $N^*$. We have
$\pi=\sum \pi e_i\, x^i$. Let $b_i=\pi e_i$ if $e_i$ is even, and
$b_i=x^i$ if $e_i$ is odd. The element $b_1\cdots b_n \in S^n(\Pi
N \oplus N^*)$ is annihilated by $\pi$. We denote by
\begin{equation}\label{eq:berezinian}
D(e_1,\dots,e_n) =b_1\cdots b_n \mod \pi S(\Pi N \oplus N^*)
\end{equation}
its class in $\Ber (N)$. We recall  the following fundamental fact
from  Manin \cite{Manin1997}*{ch. 3.4.7}:
\begin{lemma}\label{lem:berezinian}
The $\K$-module  $\Ber(N)$ is free of rank $1$, with basis
$D(e_1,\dots,e_n)$.
\end{lemma}
The module $\Ber(N)$ is a generalization of the top part of the
exterior algebra in the free ungraded situation. Indeed, when
$\K=\K_\0$ and $N=N_\0$ is free with a finite basis $(e_i)_{1\leq i\leq
n}$
consisting entirely of even elements, then $\Ber(N)$ is the top
degree part of the exterior algebra of $N$.
\begin{remark}\label{rem:berezinian}
Manin's definition of the Berezinian module is in fact $\Pi^p
\Ber(N)$, where $(p,q)$ is the rank of $N$ --- see
\cite{Manin1997}*{ch. 3.4.7}. In our presentation, there is no
advantage in doing this change of parity.
\end{remark}

\medskip
If $A$ is an associative algebra with unit, we denote by
$A^\times$ the multiplicative  group of invertible elements. We
denote by $GL(N)=\End(N)_\0^\times$ the group of even isomorphisms
of $N$. The natural action of $g\in \GL(N)$ in $\Ber(N)$ is by an
even invertible scalar $\Ber(g)\in \K_\0^\times$  (or $\Ber_N(g)$
if we need to  specify) called the \emph{Berezinian of $g$}.

Similarly, the natural  Lie action of $X \in \End(N)$ in $\Ber(N)$
is denoted by $\str(X)$ (or $\str_N(X)$), and $\str \in
\End(N)^*_\0$ is a $1$-dimensional representation of $\End(N)$
(considered as a Lie superalgebra) called the \emph{supertrace}.
In terms of a basis, we have (see \cite{Manin1997}*{ch. 3.3})
\begin{equation}\label{eq:berezinianmat}
\str(X)=\sum (-1)^{p(e_i)(p(e_i)+p(X))} x^i(X e_i).
\end{equation}

\medskip
Berezinian modules have   functorial properties which  are obvious
on the level of the defining complexes: $\Ber(N)\simeq \Ber((\Pi
N)^*)$, $\Ber(N\oplus N')\simeq \Ber(N)\otimes \Ber(N')$.
Moreover, if $\K'$ is a commutative $\K$-superalgebra, we have a
canonical isomorphism
\begin{equation}\label{eq:extber}
\Ber_{\K'}(N\otimes_\K \K')\simeq \Ber_\K(N)\otimes_\K \K'.
\end{equation}
It is less obvious that there is a natural duality between
$\Ber(N^*)$ and $\Ber(N)$:
\begin{lemma}\label{lem:dual}
There exists an unique $\K$-bilinear map
$<\cdot,\cdot>:\Ber(N^*)\times \Ber(N)\to \K$ such that, for any
basis $(e_i)_{1\leq i\leq n}$ of $N$, denoting by  $(x^i)_{1\leq
i\leq n}$ the dual basis of $N^*$, we have
\begin{equation}\label{eq:duality}
<D(x^n,\dots,x^1),D(e_1,\dots,e_n)>=1.
\end{equation}
\end{lemma}
\begin{remark}
Note that we reversed the order of the basis of $N^*$ to avoid
signs. The following computational proof  of\ lemma \ref{lem:dual}
is not very informative. It would be interesting to replace it by
a more categorical argument.
\end{remark}
\begin{proof}
 The
lemma is equivalent to the following formula:
\begin{equation}\label{eq:transber}
 \Ber(g)=\Ber(g^*) \mbox{ for } g \in GL(N),
\end{equation}
where $g^* \in \End(N^*)$ is the transpose of $g$. This identity
is given in \cite{Manin1997}*{ch. 3.3.6}.
\end{proof}
In other terms, there is a canonical isomorphism
\begin{equation}\label{eq:canonical8}
\Ber(N^*)\otimes \Ber(N)\to \K.
\end{equation}
We denote by
\begin{equation}\label{eq:canonical9}
\mathbf d:=D(x^n,...,x^1)\otimes D(e_1,...,e_n)\in
\Ber(N^*)\otimes \Ber(N)
\end{equation}
the inverse image of $1$, and we call it the \emph{canonical
element}.

\medskip Let $M$ be a $\K$-module.  Since $N$ is free of finite
rank, the injection $N\subset F(M,N)$ induces isomorphisms
$F(M)\otimes N \simeq F(M,N)$ and $F(M,\End_\K(N)) \simeq
\End_{F(M)}(F(M,N))$. Let $X \in F(M,\End_\K(N))_\0$ be such that
$\delta(X )=0$. We have (see \cite{Manin1997}*{ch. 3.3} for a
similar statement):
\begin{equation}\label{eq:berc}
\Ber_{F(M,N)}(\exp(X ))=\exp(\str_{F(M,N)}(X )).
\end{equation}
Note that the two sides of this equation are well defined elements
of $F(M)$. For the left hand side, this is because $\exp(X )\in
\GL(F(M,N))$. For the right hand side, this is because
$\delta(\str(X ))=0$.

Consider in addition  $r=\sum_{n\geq 0} r_n t^n \in \K_\0[[t]]$
such that $r_0=1$. Since $\delta(r(X))=1$, the element $r(X)$ is
invertible and has a Berezinian  $\Ber(r(X)) \in F(M)$. Define
$w=\log(r) \in t\K_\0[[t]]$ by the formula:
\begin{equation}\label{eq:berh}
w=\log(1+( r-1)) =(r-1) -\frac{(r-1)^2}{2}+\cdots=r_1 t + (r_2
-\frac{r_1^2}{2}) t^2+\cdots
\end{equation}
By functional calculus, we have $r(X)=\exp(w(X))$ and, using
(\ref{eq:berc}), we find
\begin{equation}\label{eq:berk}
\Ber(r(X))=\exp(\str(w(X)).
\end{equation}

\subsection{Dualizing modules.}\label{sec:canonical}
Let $y^1,..., y^n$  be  formal  even or odd variables. We denote
by $\K[[y^1,\dots,y^n]]$ the superalgebra of formal series, that
is the commutative superalgebra of formal sums with coefficients
in $\K$ of monomials $(y^1)^{m_1}\cdots (y^n)^{m_n}$, with $m_i
\in \N$, and $m_i\leq 1$ if $y^i$ is odd. If $p$ is the number of
even variables and $q$ the number of odd variables, we say that
the dimension of $\K[[y^1,\dots,y^n]]$ is $(p,q)$. We denote by
$\delta \in \Hom(\K[[y^1,\dots,y^n]], \K)$ the algebra
homomorphism which is zero on the $y^i$.
\begin{remark}
\label{rem:augmented} Let $\K$ be a field. Then the algebra
$\K[[y^1,\dots,y^n]]$ is local, and $\delta$ is its  unique
algebra homomorphism with values in $\K$. As we prefer not to make
this assumption on $\K$, we specify $\delta$ and consider
$\K[[y^1,\dots,y^n]]$ as an augmented algebra.
\end{remark}

{\bf Augmented algebras.} Recall that an associative
$\K$-superalgebra $A$ with unit is called \emph{augmented} if it
is provided with a $\K$-algebra morphism $\delta : A \to \K$, the
augmentation. We use the notation $A^+=\ker \delta\subset A$. The
map $\delta$ induces an isomorphism $A/A^+ \simeq \K$ and provides
$\K$ with an $A$-algebra structure.

We provide $A$ with the $A^+$-adic topology: a linear endomorphism
$X\in \End(A)$ is continuous if for all $n\in \N$, there exists
$m\in \N$ such that $X((A^+)^m)\subset (A^+)^n$. The $A$-module of
continuous derivations of $A$ is denoted by $\CT(A)$. A linear
form $w\in A^*$ is continuous if there exists $n\in \N$ such that
$w|_{(A^+)^n}=0$. We denote by $A_*\subset A^*$ the continuous
dual (it is also called the space of \emph{distributions with
support $ \{\delta\}$}).

If $\CV$ is an $A$-module, we consider the $\K$-module
$V:=\CV/A^+\CV=\CV\otimes_A \K$, and we denote by $\bar
\delta=\id_\CV\otimes_A \delta$ the quotient map $\CV \to V$. We
denote by $\a$ the $\K$-module of $\delta$-derivations $A \to \K$.

\begin{example} {\bf (Basic example).}
Let $A=\K[[y^1,\dots,y^n]]$.
We use the notation
$\frac{\partial}{\partial y^i}$ for the continuous derivation of
$\K[[y^1,\dots,y^n]]$ such that $\frac{\partial y^j}{\partial
y^i}=\delta_i^j$. These derivation form a basis of the
$\K[[y^1,\dots,y^n]]$-module $\CT(\K[[y^1,\dots,y^n]])$. The
elements $e_i:=(-1)^{p(y^i)}\delta \circ\frac{\partial}{\partial
y^i}$ provide a basis of $\a$. The images of $y^i$ in
$A^+/(A^+)^2$ form the dual basis  of $\a^*$, and the module
$\bigoplus_{i=1}^n \K y^i$ is a supplementary factor for $(A^+)^2$ in
$A^+$.
\end{example}

\begin{definition} \label{def:smooth}
We say that an augmented $\K$-superalgebra $A$ is formally smooth
of dimension $(p,q)$ if it is isomorphic to an augmented algebra
$\K[[y^1,\dots,y^n]]$ of dimension $(p,q)$.
\end{definition}
Thus, if $A$ is formally smooth, the choice of an isomorphism
$\theta^*: A \to \K[[y^1,\dots,y^n]]$ is equivalent to the choice
of the corresponding elements $x^i:={\theta^*}^{-1}(y^i) \in A^+$.
We write $A=\K[[x^1,\dots,x^n]]$, and call the family
$(x^i)_{1\leq i\leq n}$ ``a system of local coordinates at the
point $\{\delta\}$''.

\medskip
{\bf The dualizing module of a formally smooth augmented algebra.}
For the rest of this section, we consider a formally smooth
augmented $\K$-superalgebra $A$   of dimension $(p,q)$. Thus $\a$
is free of rank $(p,q)$ over $\K$ and  $\CT(A)$ is free of rank
$(p,q)$ over $A$. The space $A_*$ is a cocommutative coalgebra
with counit $\delta$, and $A$ can be identified to the dual of
$A_*$.

The dual $\Hom_A(\CT(A),A)$ is denoted by $\Omega^1_{ev}(A)$: it
is the space of \emph{even $1$-differential forms}. We denote by
$\Omega^1(A)$ the space $\Pi \Omega^1_{ev}$: it is the space of
\emph{$1$-differential forms}. For $f \in A$, we denote by
$d_{ev}f$ the element of $\Omega^1_{ev}(A)$ such that $\CT(A)\ni
\zeta\mapsto \pm \zeta(f)$, and $df:=\pi d_{ev} f \in
\Omega^1(A)$.

Following \cite{MPV1990}, we define the \emph{dualizing module}
$K(A)$ by the formula
\begin{equation}\label{eq:canonique}
K(A)=\Ber_{A}(\Omega_{ev}^1(A)).
\end{equation}
It is a free $A$-module of rank $1$.

Let $(x^i)_{1\leq i\leq n}$ be a choice of local coordinates at
the point $\{\delta\}$. The family  $(d_{ev} x^i)$ is a basis of
$\Omega^1_{ev}(A)$ with dual basis of $\CT(A)$ the vector fields
$\frac{\partial}{\partial x^i}$. Recall that  --- see formula
(\ref{eq:berezinian})---  the class of $b_1\cdots b_n$, with
$b_i=dx^i$  if $x^i$ is even, or $b_i=\frac{\partial}{\partial
x^i}$ if $x^i$ is odd, is a basis of $K(A)$ which is denoted by
$D(d_{ev}x^1,\dots,d_{ev}x^n)$.

Thus, the $\K$-module $K(A)/A^+ K(A)$ is free of rank $1$.

\begin{remark}
Since $\bar\delta:\Omega_{ev}^1(A)/A^+\Omega_{ev}^1(A) \to \a^*$
is an isomorphism,  we obtain from (\ref{eq:extber}) an
isomorphism
\begin{equation}\label{eq:canonique4}
\bar\delta : K(A)/A^+ K(A) \to \Ber_\K(\a^*).
\end{equation}
In coordinates, for $\phi \in A$, we have
$\bar\delta(D(d_{ev}x^1,\dots,d_{ev}x^n)\phi)=D(x^1,\dots,x^n)\delta(\phi)$.
\end{remark}

{\bf Tensorial modules and divergence.} Recall that $\CT(A)$ is a
Lie $\K$-subalgebra of $\End_\K(A)$.
\begin{definition}\label{def:TAmodule}
We say that an $A$-module $\CV$ endowed with a $\CT(A)$-module
structure is a $\CT(A)-A$-module if we have
\begin{equation}\label{eq:tensor}
\zeta(f v)=\zeta(f) v \pm f \zeta(v) \mbox{ for all } f\in A, v\in
\CV, \zeta \in \CT(A).
\end{equation}
The action of $\zeta $ is called the Lie derivative with respect
to $\zeta$ in $\CV$, and is traditionally denoted  by
$\CL(\zeta)$.
\end{definition}
If $\CV$ is a $ \CT(A)-A$-module, so are the module $\Pi \CV$ and
the dual $\CV^*=\Hom_A(\CV,A)$, each one with the natural actions.
If $\CV$ and $\CW$ are $ \CT(A)-A$-modules, so are $\CV\otimes_A
\CW$ and $\Hom_A(\CV,\CW)$, with the natural actions.

For instance, $A$ with the Lie derivative $\CL(\zeta)=\zeta$ is a
$ \CT(A)-A$-module.
\begin{example}
Another important example for this text  is $\CT(A)$, with the Lie
derivative $\CL(\zeta)=\ad(\zeta)$.
\end{example}
\noindent It follows that the tensor spaces, obtained from
$\Omega_{ev}^1(A)$, $\Omega^1(A)\equiv\Pi \Omega_{ev}^1(A)$,
$\CT(A)$ and $\Pi\CT(A)$ by taking tensor products over $A$  and
$\CT(A)$-invariant subquotients, are naturally $
\CT(A)-A$-modules. We  call these $\CT(A)-A$-modules
\emph{tensorial modules}.
\begin{remark}
In particular, since $\pi\in \Hom_A(\Omega_{ev}^1(A),\Omega^1(A))$
is annihilated by the Lie derivatives, $K(A)$ inherits a Lie
derivative and is a tensorial module.
\end{remark}
We shall need the following relation (see \cite{Kosmann2000} for a
similar statement):
\begin{equation}\label{eq:canonique2}
\CL(f\zeta)(D)=f \CL( \zeta)(D)\pm \zeta(f)D, \mbox{  for  } f \in
A, \zeta\in \CT(A), D\in K(A) .
\end{equation}

\medskip
 Let $D \in K(A)$ be a basis.  The \emph{divergence
$\diverg_D(\zeta)\in A$  of $\zeta\in \CT(A)$ with respect to  the
basis $D$} is defined by
\begin{equation}\label{eq:divergence}
\CL(\zeta)(D)=\diverg_D(\zeta) D.
\end{equation}

\medskip
{\bf Invariance by translation.} Let $M$ be a free $\K$-module of
finite rank $(p,q)$. We specialize to the case $A=F(M)$. Let
$(e_i)_{1\leq i \leq n}$ be a  totally ordered basis of $M$, and
$(x^i)_{1\leq i \leq n}$ the corresponding dual basis of $M^*$.
This provides an identification $F(M)=\K[[x^1,\dots,x^n]]$, and
$F(M)$ is formally smooth with local coordinates $(x^i)_{1\leq i\leq n}$.

We say that an element $d\in K(F(M))$  is \emph{invariant by
translations} if $\CL(\zeta_m) d=0$ for all $m\in M$ ($\zeta_m$ is
the derivative in the direction $m$ defined in (\ref{eq:5})).

\begin{lemma}\label{lem:haar}
The $\K$-submodule of $K(F(M))$ of elements which are invariant by
translations is free of rank $1$. A basis of this module is also
a basis of $K(F(M))$ considered as an $F(M)$-module.
\end{lemma}
\begin{proof}
In coordinates, $D(d_{ev}x^1,\dots,d_{ev}x^n)$ is such a basis.
\end{proof}
For  a basis $d$ invariant by translation, we write simply
$\diverg$ for $\diverg_d$. Let $\sum \zeta^i
\frac{\partial}{\partial x^i}$ be a derivation of
$\K[[x^1,\dots,x^n]]$. The following well known formula follows
from invariance by translations and formula (\ref{eq:canonique2}):
\begin{equation}\label{eq:divergence2}
\diverg(\sum \zeta^i \frac{\partial}{\partial x^i})=\sum \pm
\frac{\partial\zeta^i}{\partial x^i}.
\end{equation}

\medskip
{\bf Isomorphisms of augmented algebras}. Let $A$ and $B$ be two
formally smooth augmented superalgebras and $\theta^*: B \to A$ an
isomorphism of augmented algebras. By transposition, $\theta^*$
induces an isomorphism of Lie $\K$-superalgebras $\theta_* :
\CT(A) \to \CT(B)$, and of $\K$-modules $\theta^* :
\Omega_{ev}^1(B) \to \Omega_{ev}^1(A)$. We   denote also  by
$\theta^* : \CT(B) \to \CT(A)$ the inverse of the map  $\theta_*$.
This induces $\K$-isomorphisms $\theta^*$ between various tensor
spaces  related to $B$ and $A$, and in particular a bijection
\begin{equation}\label{eq:inverse}
\theta^*:K(B)\to K(A)
\end{equation}
which satisfy the relations
\begin{equation}\label{eq:inverse4}
\theta^*(f \omega)=\theta^*(f)\theta^*(\omega)\mbox{  for } f\in B
\mbox{  and  } \omega\in K(B),
\end{equation}
\begin{equation}\label{eq:inverse6}
\theta^*(\zeta \omega)=\theta^*(\zeta)\theta^*(\omega)\mbox{  for
} \zeta \in \CT(B) \mbox{  and  } \omega\in K(B).
\end{equation}

\subsection{Dualizing modules and duality.}\label{sec:odd}
Let $A$ be a formally smooth augmented superalgebra. Among the
many possible tensorial modules, the dualizing module $K(A)$ is
specially important because of its  role in integration and
duality. This is explained  in the case of supermanifolds
(algebraic, complex analytic, smooth) in
\citelist{\cite{Manin1997}\cite{MPV1990}\cite{Voronov1991}}. In
our formal context, see the work of Chemla
\citelist{\cite{Chemla1993}\cite{Chemla1994}\cite{Chemla1994b}}.
To study Gorelik elements, we need only the simple particular case
which we describe below.

\medskip
Assume that $A$ has  dimension $(0,q)$. Let $x^1,  \dots ,x^q$ be
odd coordinates. We have $A=\K[x^1,\dots,x^q]$ ---in this case
formal series are polynomials. As a $\K$-module, $A$ is free of
finite rank $2^q$. It implies that the dual $A^*$ is equal to the
continuous dual $A_*$. We have $K(A)=S^q(\CT(A))$. Recall that the
vector fields $ \frac{\partial}{\partial
x^{1}},\dots,\frac{\partial}{\partial x^{q}} $ are a basis of
$\CT(A)$. Any element of $K(A)$ can be written in an unique manner
as $ \frac{\partial}{\partial x^{1}} \cdots
\frac{\partial}{\partial x^{q}}\, \phi $ with $\phi \in A$.
Berezin integral $\int{} : K(A) \to \K$ is defined by the formula
\begin{equation}\label{eq:integration}
\int \frac{\partial}{\partial x^{1}} \cdots
\frac{\partial}{\partial x^{q}}\, \phi =\delta\left(
\frac{\partial}{\partial x^{1}} \cdots \frac{\partial}{\partial
x^{q}}( \phi)\right ).
\end{equation}
The main point, due to Berezin,  is that Berezin integral  does
not depend on the choice of the coordinates  $x^1,  \dots ,x^q$.
Berezin integral provides an isomorphism of $ \CT(A)-A$-modules
$\psi_A: K(A)\to A_*$ by the formula
\begin{equation}\label{eq:integration5}
\psi_A(\omega)(\phi)=\int \omega \phi \ \ \mbox{     for }
\omega\in K(A), \phi \in A.
\end{equation}

\medskip
The choice of coordinates identifies  $A$ and $S(\a)^*$, and $A_*$
and $S(\a)$. We write $\psi_A$ using this identification, and
the elements $e_i\in\a$ such that $x^j(e_i)=\delta^j_i$ ---
note that $e_i=-\delta\circ\frac{\partial}{\partial x^{i}}$.

\begin{lemma}
\label{lem:integration6} We have
$$
\psi_A\left(  \frac{\partial}{\partial x^{1}}
       \cdots \frac{\partial}{\partial x^{q}}
           \, \phi
\right) = (-1)^q e_1\dots e_q \phi .
$$
\end{lemma}
\begin{proof}
Since $\psi_A$ is a morphism of $A$ modules, it is sufficient to
prove the lemma for $\phi=1$. We have $\psi_A(
\frac{\partial}{\partial x^{1}} \cdots \frac{\partial}{\partial
x^{q}} )(\phi')=0$ if $\phi'\in \K[x^1,\dots,x^q]$ is of degree
smaller then $q$, and
 $\psi_A( \frac{\partial}{\partial x^{1}} \cdots
\frac{\partial}{\partial x^{q}} )(x^q\cdots x^1)=1$.  This proves
the lemma.
\end{proof}

\medskip
Let $B$ another augmented algebra and $\theta^*: B \to A$ and
isomorphism of augmented algebras. It is the transpose of a
coalgebra isomorphism $\theta : A_* \to B_*$. Recall from
(\ref{eq:inverse}) the isomorphism $\theta^*:K(B) \to K(A)$.
\begin{lemma}\label{lem:dualite2}
We have $\psi_B=\theta\circ\psi_A\circ \theta^*$.
\end{lemma}
\begin{proof}
For $\omega \in K(B)$, we have $\int \theta^*(\omega)=\int
\omega$. The lemma follows from (\ref{eq:inverse4}).
\end{proof}

\section{Dualizing modules  on formal Lie supergroups.}\label{sec:supergroups}
Let $\g=\g_\0\oplus \g_\1$ be a Lie $\K$-superalgebra, and $U(\g)$ be
its enveloping algebra. We denote by $j: \g \to U(\g)$ the
canonical map. Since $\K_\0$ contains $\Q$, $j$ is injective (cf.
\citelist{\cite{Cohn1963} \cite{Petracci2003}}).

\subsection{Enveloping algebras as coalgebras.}
\label{sec:newsection}
 We
use the notation $U^+(\g)$ for the ideal generated by $\g$.

Recall  that $U(\g)$ is a cocommutative Hopf superalgebra. For  $u
\in U(\g)$, we write
\begin{equation}\label{eq:coprod}
\Delta_{U(\g)}(u)=\sum u_i\otimes u'_i
\end{equation}
for its coproduct, and we denote by $S_{U(\g)} : U(\g)\to U(\g)$
the antipode.  We use the notation $F(G)$ for the commutative
superalgebra $U(\g)^*$  (dual of the commutative
coalgebra $U(\g)$) provided with the augmentation $\delta: f \to f(1)$.
This can be considered as the {\it definition of a formal Lie
supergroup $G$}:  $F(G)$ is the space of formal functions defined
in a formal neighborhood of $1\in G$, and $\delta$ is the
evaluation at $1$.

For a $\K$-module $N$, we define in a similar way
$F(G,N)=\Hom(U(\g),N)$, considered as the $N$-valued functions on
$G$. As in subsection \ref{sub:1}, $F(G,N)$ is an $F(G)$-module,
and, if $\mu :N_1\otimes N_2 \to N$ is a map of $\K$-modules,
there is a corresponding map $\tilde\mu: F(G,N_1)\otimes_{F(M)}
F(G,N_2)\to F(G,N)$. In particular, if $N$ is a $K$-superalgebra,
$F(G,N)$ is an $F(G)$-superalgebra.

\medskip
If $f_1$ and $f_2$ are in $\End(U(\g))$ we use the notation
$f_1*f_2$ for the product in $ F(G,U(\g))=\End(U(\g))$. Formulas
(\ref{eq:6})  generalize in two versions, one for {\it left and} one for
{\it right action}. Here are the details. Let $\phi \in \End(U(\g))$. We
use the notation $c^\lambda_\phi=\phi *\id_{U(\g)}$ and we denote
by $\mathbf \zeta^\lambda_\phi$ the transpose map acting in
$F(G,N)$. For $u\in U(\g)$ and $g\in F(G,N)$,  we have
\begin{equation}\label{eq:13}
c^\lambda_\phi(u)=\sum \phi(u_i)u'_i \mbox{ and }
\zeta^\lambda_\phi(g)=\pm g\circ c^\lambda_\phi.
\end{equation}
Similarly, we  use the notation $c^\rho_\phi= \id_{U(\g)} * \phi $
and we denote by $\zeta^\rho_\phi$ the transpose map acting in
$F(G,N)$. Thus we have
\begin{equation}\label{eq:15}
c^\rho_\phi(u)=\sum \pm u_i \phi(u'_i) \mbox{ and }
\zeta^\rho_\phi(g)=\pm g\circ c^\rho_\phi.
\end{equation}

\medskip
Note that $\phi=S_{U(\g)}*c^\rho_\phi=c^\lambda_\phi*S_{U(\g)}$,
and that $c^\rho_\phi=\id_{U(\g))}*c^\lambda_\phi*S_{U(\g)}$.

\medskip
Let $v\in U(\g)$. We denote by $l_v$ the left multiplication by
$v$ in $U(\g)$. Identifying $v$  with the constant function on $G$
with value $v$ we get $c^\lambda_v=l_v$. Similarly,
$c^\rho_v(u)=\pm uv$, $\zeta^\lambda_v$ is the right invariant
differential operator whose evaluation at $1$ is equal to $v$, and
$\zeta^\rho_v$ is the left invariant differential operator whose
evaluation at $1$ is equal to $v$. For $a\in \g$, we will simply
write $l_a:=l_{j(a)}$, $\zeta^\lambda_a:=\zeta^\lambda_{j(a)}$ and
$\zeta^\rho_a:=\zeta^\rho_{j(a)}$.

\medskip
Let $f\in F(G)$. Identifying  it with a scalar valued element  of
$\End(U(\g))$, we get  $c^\rho_f=c^\lambda_f$. This provides
$U(\g)$ with a structure of $F(G)$-module. Accordingly, for $u\in
U(\g)$, we also write $fu$ for $c_f(u)$.

\medskip
An element $\alpha \in F(G,\g)$ is called a \emph{vector field on
$G$}.
\begin{remark}
\label{rem:asterisco} Both $\zeta^\lambda_\alpha$ and
$\zeta^\rho_\alpha$ are elements of $\CT(F(G))$.
\end{remark}
Let $\beta \in F(G,\g)$ be a vector field. We have
\begin{equation}\label{eq:17}
[\zeta^\rho_\alpha, \zeta^\rho_\beta]=\zeta^\rho_\gamma  \mbox{
with } \gamma= \zeta^\rho_\alpha(\beta)-\pm
\zeta^\rho_\beta(\alpha)+[\alpha, \beta],
\end{equation}
\begin{equation}\label{eq:17a}
[\zeta^\lambda_\alpha,
\zeta^\lambda_\beta]=\zeta^\lambda_{\gamma'}
 \mbox{ with }
 \gamma'=
\zeta^\lambda_\alpha(\beta)-\pm
\zeta^\lambda_\beta(\alpha)-[\alpha, \beta],
\end{equation}
\begin{equation}\label{eq:17b}
[\zeta^\lambda_\alpha, \zeta^\rho_\beta]=\zeta^\rho_{\gamma_1}-\pm
\zeta^\lambda_{\gamma_2},
 \mbox{ with }
\gamma_1=  \zeta^\lambda_\alpha(\beta) \text{  and  } \gamma_2=
\zeta^\rho_\beta(\alpha).
\end{equation}
These formulas include as particular cases formulas
(\ref{eq:11}) when $\g$ is commutative, and the formulas
$[\zeta^\rho_a, \zeta^\rho_b]=\zeta^\rho_{[a,b]}$,
$[\zeta^\lambda_a, \zeta^\lambda_b]=-\zeta^\lambda_{[a,b]}$ and
$[\zeta^\lambda_a,\zeta^\rho_b]=0$ for $a\in \g$ and $b\in \g$. We
have similar relations for the coderivations $c^\rho_\alpha$ and
$c^\lambda_\alpha$ of $U(\g)$. We note the following classical
proposition.

\begin{proposition}\label{prop:coder}
The map  $\alpha \to c^\rho_\alpha$ is an isomorphism from
$F(G,\g)$ to the space of coderivations of $U(\g)$.
\end{proposition}
\begin{proof}
Since   $\K$ contains $\Q$, $\g$ is the subspace of primitive
elements of $U(\g)$.  As in \cite{Radford1986}*{theorem~1}, this
implies the proposition.
\end{proof}

\subsection{The exponential map of a formal Lie supergroup.}
Recall from \cite{Cohn1963}, \cite{Petracci2003} that the
symmetrization $\beta: S(\g)\to U(\g)$ is  a coalgebra isomorphism
such that $\beta (1)=1$  and $\beta(a)=j(a)$    for $a\in \g$.
\begin{remark}
The transpose map $\beta^*: F(G) \to F(\g)$ is an isomorphism of
augmented algebras.
\end{remark}
We denote also by  $\beta^*: F(G,N) \to F(\g,N)$ the map defined
by $\beta^*(g)=g \circ \beta$ for $g\in F(G,N)$.

We consider $\beta $ as a formal diffeomorphism   from a
neighborhood of $0\in \g$ to a neighborhood of $1\in G$. This is
the formal analog of the usual exponential map. It can also be
written as an exponential in our formal setting: extend $j$ to an
element of the algebra $F(\g,U(\g))$ by assigning the value zero
on $S^n(\g)$ for $n\neq 1$. The sum $\exp(j)=\sum_{n\geq 0}
\frac{j^n}{n!}$ is   locally finite and $\beta=\exp(j)$.

\subsection{The generic point of a Lie
superalgebra.} \label{sec:generic} In this subsection, we collect
some results on $F(\g)$ needed later.  We write simply $x=x_\g \in
F(\g,\g)$ the generic point of $\g$.

We consider $\ad x \in \End_{F(\g)}(F(\g,\g)) $. For
$p=\sum_{n\geq 0} p_n t^n \in \K_\0[[t]]$, we can define  $p(\ad
x)\in \End_{F(\g)}(F(\g,\g)) $ by a locally finite sum as in
(\ref{eq:function}). For any vector field $\alpha \in F(\g,\g)$,
we obtain a new vector field $p(\ad x)(\alpha)$. As in
\cite{Petracci2003}, the vector fields
\begin{equation}\label{eq:rxa}
 p(\ad x)(a)
=\sum_{n\geq 0}
 p_n  \ad^n (x )(a) \mbox{  with  } a\in \g,
\end{equation}
their corresponding coderivations $c_{p(\ad x)(a)}$ of $S(\g)$,
and the transpose operators  $\zeta_{p(\ad x)(a)}$ of $F(\g,N)$
($N$ is a $\K$-module), play a special role. We recall how to
compute $p(\ad x)(a)$. Let $b_1$, $b_2$, \dots, $b_n$ in $\g$. We
have (see \cite{Petracci2003})
\begin{equation}\label{eq:radx}
        p(\ad x)(a) : b_1\cdots b_n
\mapsto p_n
        \sum_{s\in S_n}
        \pm
        \ad b_{s(1)}\cdots \ad b_{s(n)} (a).
\end{equation}

\medskip
In the rest of this  section, we suppose that $\g$ is free of
finite rank. The Lie superalgebra  $F(\g,\g)=\g\otimes F(\g)$ is
free of finite rank over $F(\g)$. Thus,   the divergence
$\diverg(\zeta_{p(\ad x) a })\in F(\g)$  of the derivation
associated to the vector field $p(\ad x)(a)$ is defined in
(\ref{eq:divergence2}). Similarly, for $w\in \K_\0[[t]]$, the
function $\str(w(\ad x)) \in F(\g)$ is defined. We state some
results relating these objects.
\begin{proposition}\label{prop:div}
\begin{equation}
  \diverg(\zeta_{p(\ad x) a })
= -\str\left( \frac{ p(\ad x)-p(0) }{\ad x} \ad a
       \right)
.
\end{equation}
\end{proposition}
\begin{proof}
It is sufficient to prove it for $p(t)=t^n$. We pick a basis
$(e_i)_{i\in I}$ of $\g$ and the corresponding dual basis $(x^i)_{i\in I}$
of
$\g^*$. Let $\zeta_{e_i}$ the derivative in the direction $e_i$.
We have
\begin{align}
    \zeta_{e_i}\big(\ad^n( x)(a)\big)
&=\sum_{p+q=n-1}\ad^p( x)\ad(e_i)\ad^q( x)(a)\\
 &= \sum_{p+q=n-1}\ad^p( x) \left([e_i,\ad^q( x)(a)] \right)\notag\\
 &=-\sum_{p+q=n-1}\ad^p( x)  \ad(\ad^q( x)(a)) (e_i). \notag
\end{align}
Applying  formula (\ref{eq:divergence2}), we get
\begin{equation}
\diverg\left(\zeta_{\ad^n( x)( a)}\right)
=-\str\left(\sum_{p+q=n-1}\ad^p(x) \ad(\ad^q( x)(a))\right ).
\end{equation}
Using the notations $X:=\ad x$, $A:=\ad a$ and the relation
$\ad(\ad^n( x)(a))= \ad^n (X)(A)$, we obtain
\begin{equation}
\diverg\left(
\zeta_{\ad^n( x)( a)}
\right)=-\str\left(\sum_{p+q=n-1}X^p
\ad^q (X)(A)\right).
\end{equation}
Let $l$ and $r$ be the left and right multiplication by $X$ in
$\End_{F(\g)}(F(\g,\g))$. Since $\ad X=l-r$, we obtain
\begin{equation}\label{eq:divergenced}
\diverg\left(
\zeta_{\ad^n( x)( a)}
\right)=-\str\left(\sum_{p+q=n-1}l^p
(l-r)^q(A)\right).
\end{equation}
For all $Y \in \End_{F(\g)}(F(\g,\g))$ we have
$\str((l-r)Y)=\str(XY-YX)=0$.  We obtain
\begin{equation}\label{eq:divergencee}
\diverg\left(
\zeta_{\ad^n( x)( a)}
\right)=-\str(l^{n-1}(A))=-\str\left(
\ad^{n-1} (x)\ad a\right),
\end{equation}
which is what we want.
\end{proof}

\begin{proposition}
\label{prop:diffbera}
\begin{equation}\label{eq:diffbera}
  \zeta_{p(\ad x)(a)}
  \left( \str\left(w(\ad x)\right)
  \right)
 = \str\left(p(0) w'(\ad x)\ad a\right).
\end{equation}
\end{proposition}
\begin{proof}
It is sufficient to prove it for $w=t^n$. Using the same notations
and arguments  in the proof of proposition \ref{prop:div} we have
\begin{align}
     \zeta_{p(\ad x)(a)}
     \left( \str(\ad^n (x) )
     \right)
 & =\sum_{p+q=n-1}\str\left(\ad^p  (x)   \ad( p(\ad x)(a) )  \ad^q ( x) \right)\\
 &= \sum_{p+q=n-1}\str\left(X^p p(\ad X) (A)X^q\right)\\
 &=\sum_{p+q=n-1}\str\left(l^p r^q p(l-r)(A)\right)
 .
  \label{eq:divergencedd}
\end{align}
Since $p(l-r)-p(0)$ is divisible by $l-r$, we obtain
\begin{align}
     \zeta_{p(\ad x)(a)}
     \left( \str(\ad^n  x )
     \right)
 &=\sum_{p+q=n-1}\str\left(p(0) l^p r^q (A)\right)\\
 & =\sum_{p+q=n-1}\str\left(p(0) X^{p+q}A)=\str(p(0) n X^{n-1}A\right) ,
\end{align}
which is what we want.
\end{proof}

Let $r(t)=1+r_1 t+ \cdots \in \K_\0[[t]]$. Since $r(\ad x)$ is an
even invertible endomorphism of $\g\otimes F(\g)$, the Berezinian
$\Ber r(\ad x)\in F(\g)$ is defined.  Recall  from (\ref{eq:berh})
the definition of $w=\log(r)$. From (\ref{eq:berk}) we obtain
\begin{equation}\label{eq:berradx}
\Ber r(\ad x) =\exp\left(\str(w(\ad x))\right).
\end{equation}
\begin{proposition}\label{prop:diffber}
\begin{equation}
\label{eq:diffber}
\zeta_{p(\ad x)(a)}
\left( \Ber r(\ad x) \right)=
 \str\left( p(0)w'(\ad x) \ad a\right)
 \Ber r(\ad x)
.
\end{equation}
\end{proposition}
\begin{proof}
By lemma \ref{lem:dercomp} applied with  $f(t)=\exp(t)$ and
$X=\str(w(\ad x))$, we have
\begin{equation}
\label{eq:diffberb}
 \zeta_{p(\ad x)(a)}
\left( \Ber r(\ad x) \right)
 =\zeta_{p(\ad x)(a)}
  \left( \str(w(\ad x))
  \right)
  \Ber r(\ad x)
.
\end{equation}
Then the result follows from proposition \ref{prop:diffbera}.
\end{proof}

\subsection{The Jacobian of the exponential map.}
\label{sec:jacexp}
In this subsection we suppose   that $\g$ is free of
finite rank as a $\K$-module. The algebra $F(G)$ is  formally
smooth, and  the dualizing module $K(F(G))$ is a free
$F(G)$-module of rank 1.

\begin{remark}
Let $\CV$ be a $\CT(F(G))-F(G)$-module (see definition
\ref{def:TAmodule}). By composing the Lie derivative with the
representations $a \to -\zeta^\lambda_a$ and $a \to \zeta^\rho_a$
of $\g$, we get  two commuting representations of $\g$ in $\CV$
(recall remark \ref{rem:asterisco}). They are called respectively
\emph{the left and right actions of $\g$}.
\end{remark}
\begin{proposition}\label{prop:free}
The map  $F(G,\g)\ni \alpha \to \zeta^\rho_\alpha\in \CT(F(G))$ is
an isomorphism of $F(G)$-modules. It intertwines the left actions
$-\zeta^\lambda_a$ and $\CL(-\zeta^\lambda_a)$ of  $a \in \g$.
\end{proposition}
\begin{proof}
The first assertion follows from proposition \ref{prop:coder}. By
definition of the Lie derivative in $\CT(F(G))$ we have
$\CL(\zeta^\lambda_a)(\zeta^\rho_\alpha)=[\zeta^\lambda_a,\zeta^\rho_\alpha]$.
Since $\zeta^\rho_\alpha(a)=0$, it follows from formula
(\ref{eq:17b}) that we have
$[\zeta^\lambda_a,\zeta^\rho_\alpha]=\zeta^\rho_\gamma$ with
$\gamma=\zeta_a^\rho(\alpha)$.
\end{proof}

\begin{remark} \label{rem:A}
By duality, we define an isomorphism $F(G,\g^*)\ni \mu \mapsto
\omega_{ev,\rho}^\mu\in \Omega^1_{ev}(F(G))$ by
$\omega_{ev,\rho}^\mu( \zeta^\rho_\alpha)=\mu(\alpha)$. It
intertwines the left actions $-\zeta^\lambda_a$ and
$\CL(-\zeta^\lambda_a)$ of $a \in \g$. In particular, it defines a
$\K$-isomorphism $\g^*\ni f \to \omega_{ev,\rho}^f$  into the
space of left invariant differential forms.
\end{remark}

\medskip
Let $(e_i)_{1\leq i \leq n}$ be a totally ordered basis of $\g$,
and   $(x^i)_{1\leq i \leq n}$ its dual basis.  The derivations
$\zeta^\rho_i:=\zeta^\rho_{e_i}$ form a basis of $\CT(F(G))$. Let
$\omega_{ev,\rho}^i \in \Omega_{ev}^1(F(G))$ be the elements of
the dual basis. We denote by
$D^\rho=D(\omega_{ev,\rho}^1,\dots,\omega_{ev,\rho}^n)$ the
corresponding basis of $K(F(G))$. It is invariant by the left
action of $\g$ in $K(F(G))$, and (recall (\ref{eq:canonique4}))
$\bar\delta D^\rho= D(x^1,\dots,x^n) \in \Ber(\g^*)$. We obtain
the following result.
\begin{proposition}
\label{prop:haar} The $\K$-module of elements of $K(F(G))$
invariant by the left action of $\g$ is free of rank $1$. A  basis
of this module is also  a basis of $K(F(G))$ as an $F(G)$-module.
\end{proposition}
Similarly,  $d:=D(d_{ev} x^1,\dots,d_{ev}x^n)$
is a basis of
$K(F(\g))$ invariant by translations (see lemma \ref{lem:haar}), and
$\bar\delta
d=D(x^1,\dots,x^n) \in \Ber(\g^*)$. The symmetrization $\beta$
induces an isomorphism $\beta^*:K(F(G)) \to K(F(\g))$, and the
element $\beta^*(D^\rho)$ is a basis of $K(F(\g))$. There exists
an unique $J^\rho \in F(\g)_\0$
\begin{equation}\label{eq:jac}
 \beta^*(D^\rho)=J^\rho d.
\end{equation}
Since $J^\rho$ does not depend on the choice of the basis $(e_i)_{1\leq
i\leq n}$,
there is no ambiguity in the definition of $J^\rho$. The function
$J^\rho$ is called the \emph{Jacobian of the exponential map in
the left invariant frame}. The following formula is well known in
the case of Lie groups.
\begin{theorem}
\label{theo:jac} Let $x$ be the generic point of $\g$. We have
\begin{equation}\label{eq:jaca1}
  J^\rho = \Ber\left( \frac{1-e^{-\ad x}}{\ad x} \right).
\end{equation}
\end{theorem}
\begin{proof}
We give two proofs. Both  use the notations
$r(t)=\frac{1-e^{-t}}{t}$, $p(t)=\frac{t}{e^t-1}$ and
$w(t)=\log(r(t))$.

\emph{First proof}. As in remark \ref{rem:A}, replacing $G$ with
$\g$, we identify $\Omega_{ev}^1(F(\g))$ with $F(\g, \g^*)$. We
recall the  formula given in \cite{Petracci2003} for
$\beta^*(\omega_{ev,\rho}^i)\in F(\g, \g^*)$. We consider $r(\ad
x)\in F(\g, \End(\g))=\End_{F(\g)}(F(\g,\g))$ as an element of
$\GL(F(\g, \g))$. Its transpose $r(\ad x)^*$ belongs to $\GL(F(\g,
\g^*))$ and we have $\beta^*(\omega_{ev,\rho}^i)= r(\ad
x)^*(x^i)$. By definition of the Berezinian, the Jacobian
$J^{\rho}$ is equal to $\Ber (r(\ad x)^*)$. The theorem follows
from formula (\ref{eq:transber}).

\emph{Second proof}. We present this less natural proof for fun,
and because it is a simplification of the proof of theorem
\ref{thm:main} below. Let $a\in \g$. We recall the formula given
in \cite{Petracci2003} for the derivation
$\beta^*(-\zeta_a^\lambda)$ of $F(\g)$ obtained by transporting
the derivation $-\zeta_a^\lambda$ of $F(G)$.  We have
\begin{equation}\label{eq:emanuela22}
\beta^*(-\zeta_a^\lambda)= -\zeta_{p(\ad x)(a)}.
\end{equation}
Since $p(0)=1$, it is sufficient  to show that for all $a\in \g$
we have
\begin{equation}\label{eq:jaca2}
\CL( \zeta_{p(\ad x)(a)} ) \left( \Ber(r(\ad x))  d \right) =0.
\end{equation}
By   (\ref{eq:canonique2}) and (\ref{eq:divergence}) this is
equivalent to
\begin{equation}\label{eq:jaca4}
  \zeta_{p(\ad x)(a)}
  \left( \Ber(r(\ad x))
  \right)
+ \Ber\left(r(\ad x)\right)
  \diverg(\zeta_{p(\ad x)(a)})=0.
\end{equation}
By propositions \ref{prop:diffber} and \ref{prop:div}, this
reduces  to
\begin{equation}\label{eq:jaca5}
\str\left(p(0) w'(\ad x) \ad a\right)-\str\left(\frac{p(\ad
x)-p(0)}{\ad x} \ad a\right)=0 .
\end{equation}
This identity follows from the relation
\begin{equation}\label{eq:jaca6}
p(0)w'(t)-\frac{p(t)-p(0)}{t}=0,
\end{equation}
which is miraculously  true:  since  $p(0)=1$, we have
$$
p(0)w'(t)=\frac{e^{-t}}{1-e^{-t}}-\frac{1}{t}=\frac{1}{e^t-1}-\frac{1}{t}
\ \ \ \mbox{and} \ \ \
\frac{p(t)-p(0)}{t}=\frac{\frac{t}{e^t-1}-1}{t}=\frac{1}{e^t-1}-\frac{1}{t}
.
$$
\end{proof}

\section{Formal homogeneous spaces.}\label{sec:inducedformal}
In this section, we consider a Lie superalgebra $\g$  and
 a Lie subalgebra $\h\subset \g$ which is a direct factor: there
exists and we choose a $\K$-submodule $\q\subset \g$ such that
\begin{equation}\label{eq:direct}
\g=\h\oplus\q.
\end{equation}

The following lemma is a straightforward consequence of the fact
that $\beta $ is an isomorphism of coalgebras.
\begin{lemma}\label{eq:lem}
Let $\g=\p\oplus\q$ be a decomposition into a direct sum of
submodules.  The  map $S(\p)\otimes S(\q)\ni w\otimes w' \to
\beta(w)\beta(w')\in U(\g)$
 is  an
isomorphism of coalgebras.
\end{lemma}
\begin{corollary}\label{cor:betabar}
The map $\bar \beta : S(\q)\otimes U(\h)\ni w\otimes u \to
\beta(w)u\in U(\g)$  is an isomorphism of coalgebras.
\end{corollary}

\subsection{Induced and coinduced representations in exponential
coordinates.}\label{sec:induced} For $a\in \g$, we  recall that
the left multiplication  $l_a$ is a coderivation of  $U(\g)$. We
denote by $\bar l_a$ the corresponding coderivation of
$S(\q)\otimes U(\h)$. Recall the notation $\Delta_{S(\g)}(w)=\sum
w_i\otimes w'_i$ for $w\in S(\q)$.

\begin{proposition}\label{prop:formulacoder}
Let $a\in \g$. There exist $\alpha_a \in F(\q,\q)$ and
$\theta_a\in F(\q,\h)$ such that, for $w\in S(\q)$ and $u\in
U(\h)$ we have
\begin{equation}\label{eq:coderx}
\bar l_a(w\otimes u)= \sum \pm w_i\alpha_a(w'_i)\otimes u+
   \sum \pm w_i\otimes j(\theta_a(w'_i))u.
\end{equation}
\end{proposition}
\begin{proof}
We use  the convolution notation for the  product in
$F(\q,S(\q)\otimes U(\h))$.  Define  $C \in F(\q, S(\q)\otimes
U(\h))$ by
\begin{equation}
\label{eq:C} C(w)=\bar l_a(w\otimes 1) \mbox{ for } w\in S(\q).
\end{equation}
We have $\bar l_a(w\otimes u)=C(w)(1\otimes u)$,  and formula
(\ref{eq:coderx}) is equivalent to
\begin{equation}
\label{eq:C1}
 C
 = (\id_{S(\q)}\otimes \id_{U(\h)})
 * ( \alpha_a\otimes \id_{U(\h)}
   + \id_{S(\q)}\otimes j\circ \theta_a
   )
 = c_{\alpha_a}\otimes 1
 + \left(\id_{S(\q)}\otimes j\circ \theta_a \right)\circ \Delta_{S(\q)}
.
\end{equation}
To prove (\ref{eq:C1}), we argue as in \cite{Radford1986}*{theorem
1}.  We define $\gamma \in F(\q, S(\q)\otimes U(\h))$ by $\gamma=
(S\otimes 1) *C$, where $S \in F(\q,S(\q))$ is the antipode of
$S(\q)$. Since $S$ is the inverse of $\id_{S(\q)}$ for the
convolution product, we have $C=(\id_{S(\q)}\otimes 1) *\gamma$.
Because $\bar l_a$ is a coderivation, $\gamma $ takes its values
in the submodule of primitive elements of $S(\q)\otimes U(\h)$.
Since $\K$ is $\Q$-algebra, this submodule is $\q\otimes 1+
1\otimes j(\h)$. We write $\gamma =\alpha_a\otimes 1+ 1\otimes
j\circ\theta_a$ with $\alpha_a \in F(\q,\q)$ and $\theta_a\in
F(\q,\h)$.
\end{proof}
\begin{remark}
Manageable formulas for the   maps $\alpha_a \in F(\q,\q)$ and
$\theta_a\in F(\q,\h)$ are not known in general. In the case
$\h=\{0\}$, a formula for $\alpha_a$ is given in
\cite{Petracci2003} (see (\ref{eq:emanuela22}) above). In the case
of symmetric pairs, we give formulas  for   $\alpha_a$  and
$\theta_a$ in theorem \ref{theo:formula}.
\end{remark}

\medskip
Let  $(V,\chi)$ be an $\h$-module. We consider the induced
representation of $\g$, that is the  representation by left
multiplication in $U(\g)\otimes_\h V$. We define
\begin{equation}\label{eq:Ind2}
\bar \beta_\chi :S(\q)\otimes V \to U(\g)\otimes_\h V
\end{equation}
 as in corollary \ref{cor:betabar}.
For  $a\in \g$, we denote by $l_{\chi,a}$ the endomorphism of
$U(\g)\otimes_\h V$ induced by the $l_{a}$, and by $\bar
l_{\chi,a}$ the corresponding endomorphism of $S(\q)\otimes V$. It
 belongs to $F(\q,S(\q)\otimes \End(V))$,  and it follows
from (\ref{eq:C1}) that $\bar l_{\chi,a}$ is given by the formula
\begin{equation}\label{eq:coderxc}
   \bar l_{\chi,a}:
   S(\q)\ni w\mapsto
  c_{\alpha_a}(w)\otimes \id_V
+ \left( \id_{S(\q)}\otimes \chi\circ\theta_a\right)\circ
\Delta_{S(\q)}(w)(1\otimes \id_V) .
\end{equation}

\medskip
The coinduced representation is a representation in the space
\begin{equation}\label{eq:Coind1}
F(G/H,V):=\Hom_\h(U(\g), V),
\end{equation}
that is the subspace of elements $g\in F(G,V)$  such that $-\pm
g(uj(a))=a(g(u))$ for $a\in \h$ and $u\in U(\g)$. The action of
$a\in \g$ in $F(G/H,V)$  is the restriction of
$-\zeta^\lambda_{a}\in\End(F(G,V))$, we denote it by
$-\zeta_{\chi,a}$.

It follows from corollary \ref{cor:betabar} that the map
\begin{equation}\label{eq:isobetastar}
\bar \beta^*_\chi: F(G/H,V)\ni \phi \to \phi\circ \beta\in F(\q,V)
\end{equation}
is bijective.  We denote by
$
      -\bar \zeta_{\chi,a}
   := -\bar\beta^*_\chi
\circ \zeta_{\chi,a} \circ \bar{\beta^*_\chi}^{-1} $ the
corresponding endomorphism of $F(\q,V)$. We denote by $\chi(
\theta_a )$ the endomorphism $
    \phi
    \to \pm
        (\chi\circ \theta_a
\otimes \phi)\circ \Delta_{S(\q)} $ of $\End_{F(\q)}(F(\q,V))$. It
follows from (\ref{eq:coderxc}) that
\begin{equation}\label{eq:coderxd}
-\bar \zeta_{\chi,a}=-\zeta_{\alpha_a}+ \chi(\theta_a).
\end{equation}

\subsection{The exponential map  on a formal  homogeneous space}
\label{sec:32}
We recall (see formula \ref{eq:Coind1}) that $F(G/H)$ is the
algebra dual to the coalgebra $U(\g)/U(\g)\h$, i.e. the orthogonal
of $U(\g)\h$ in $F(G)$.

The restriction of $\delta$ to $F(G/H)$ is still denoted by
$\delta$. It provides $\K$ with a structure of $F(G/H)$-module.
The augmented algebra $F(G/H)$ is thought of as the algebra of
functions on the \emph{formal homogeneous space $G/H$}, and
$\delta$ is ``the evaluation at the point $1\in G/H$''.

Note that $U(\g)/U(\g)\h$ and $F(G/H)$ are the   induced and
coinduced representations by the trivial representation of $\h$ in
$\K$. We denote this representation by $\triv$, and we use the
notations $\widetilde \beta :=\bar\beta_{\triv} $ and $\widetilde
\beta^* :=\bar\beta^*_{\triv} $. In this case,
\begin{equation}\label{eq:isoq}
\widetilde \beta : S(\q) \to U(\g)/U(\g)\h \mbox{  and }
\widetilde \beta^*: F(G/H) \to F(\q) ,
\end{equation}
are respectively a  coalgebra and an algebra isomorphisms. The map
$\widetilde \beta$ is the formal analogue of  the exponential map
$\q \to G/H$.

For $a\in \g$, we still denote by $l_{a}$ the left multiplication
by $j(a)$  in the quotient $U(\g)/U(\g)\h$, i.e.
$l_{a}=l_{\triv,a}$. We denote by $\zeta_{a}$ the transpose map in
$F(G/H)$. Thus $\zeta_{a}=\zeta_{\triv,a}$ is the restriction of
$\zeta^\lambda_{a}$ to $F(G/H)\subset F(G)$.

Consider the maps in (\ref{eq:isoq}). For $a\in \g$, we denote
$\widetilde l_a:= \bar l_{\triv, a}$ and $
   -\widetilde \zeta_a
:= -\bar\zeta_{{\triv}, a} $ the corresponding coderivation of
$S(\q)$  and derivation of $F(\q)$. As a particular case of
formulas (\ref{eq:coderxc}) and (\ref{eq:coderxd}), we have
\begin{equation}
\label{eq:isoq3} \widetilde l_a=c_{\alpha_a}  \mbox{  and }
-\widetilde \zeta_a=-\zeta_{\alpha_a}
\end{equation}
where $\alpha_a \in F(\q,\q)$ is the vector field defined in
proposition \ref{prop:formulacoder}.

\subsection{Imprimitivity theorem.}\label{sec:imprimitivity}
The following definition is similar to  definition
\ref{def:TAmodule}.
\begin{definition}\label{def:Ind}
A $\g-F(G/H)$-module $\CV$ is a $\K$-module provided with a
structure of $F(G/H)$-module and a structure of $\g$ module which
satisfy the Leibnitz relation
\begin{equation}\label{eq:Leibnitz}
  a(\phi f)
=- \zeta_{a}^\lambda (\phi)f \pm \phi(af), \mbox{  for } a\in \g,
\phi\in F(G/H), f\in \CV .
\end{equation}
\end{definition}
If $\CV$ and $\CW$ are $\g-F(G/H)$-modules, so are $\CV
\otimes_{F(G/H)} \CW$ and $\Hom_{F(G/H)}(\CV,\CW)$.

\medskip
Induced and coinduced modules by a representation of $\h$ in a
$\K$-module $V$ are  naturally $\g-F(G/H)$-modules. Indeed: as
$F(G/H,V)\subset F(G,V)$ is stable under multiplication by
$F(G/H)$, then $F(G/H,V)$ is a $\g-F(G/H)$-module. We still denote
by $\delta$ the evaluation map $\delta:F(G/H,V)\to V$. We remark
that we have natural isomorphisms of $\g-F(G/H)$-modules
\begin{equation}\label{eq:isoisoa}
U(\g)/U(\g)\h \otimes_{F(G/H)}  F(G/H,V) \simeq U(\g)\otimes_\h V,
\end{equation}
and
\begin{equation}\label{eq:isoisob}
F(G/H,V) \otimes_{F(G/H)}  F(G/H,W) \simeq F(G/H,V\otimes W).
\end{equation}
The purpose of imprimitivity theorems is to characterize induced
and coinduced  modules among  the $\g-F(G/H)$-modules. We state
such a theorem for coinduced representations in the particular
simple case we need.

\medskip
Let $\CV$ be a $\g-F(G/H)$-module. We consider
$V_{\CV}:=\CV\otimes_{F(G/H)}\K$ and the corresponding map $\bar
\delta:=\id_\CV\otimes_{F(G/H)} \delta : \CV \to V_{\CV}$. The
following straightforward lemma  provides $V_{\CV}$ with a natural
$\h$-module structure.
\begin{lemma}\label{lem:quotient}
The kernel of $\bar \delta$ in $\CV$ is stable under the action of
$\h$.
\end{lemma}
  Suppose also  that $V$ is a free
$\K$-module of finite rank $(p,q)$. It follows that $F(G/H,V)$ is
a free $F(G/H)$-module of finite rank $(p,q)$. Moreover, the map
$\delta :F(G/H,V) \to V$ induces an isomorphism of $\h$-modules
$V_{F(G/H,V)} \to V$, that is we have $\bar \delta=\delta$. In
particular,
$$
F(G/H, V) \simeq F(G/H, V_{F(G/H,V)}) .
$$
We need a converse to this property, the imprimitivity theorem,
due to Blattner \cite{Blattner1969} (see also
\cite{Scheunert1979}) in a slightly different setting.

\medskip
Assume that $\CV$  is free of finite rank  as an $F(G/H)$-module.
The   $\K$-module $V_{\CV}$ is free of the same rank. Define a map
\begin{equation}
\label{eq:coinducemap} \tilde \delta: \CV \to F(G/H,V_{\CV})
\mbox{ by }\tilde \delta(v)(u)=\pm \bar \delta(uv).
\end{equation}
The facts that $\tilde \delta$ is well defined, and that it
commutes with the actions of $\g$ and $F(G/H)$ are easy.
\begin{proposition}
\label{prop:blattner} Let $\CV$ be a $\g-F(G/H)$-module which is
free of finite rank  as an $F(G/H)$-module. Then  $\tilde \delta$
is an isomorphism. In particular
$$
\CV \simeq F(G/H, V_{\CV}) .
$$
\end{proposition}
\begin{proof}
Let $(\phi_i)$ be a basis of $\CV $, and  $v_i:=\bar
\delta(\phi_i)\in V_{\CV}$. Then $(v_i)$ is a basis of $V_{\CV}$.
We have $\delta(\tilde \delta(\phi_i))=v_i$. It follows from lemma
\ref{lem:basis} that $(\tilde \delta(\phi_i))$ is a basis of
$F(G/H,V_{\CV})$.
\end{proof}

\subsection{Dualizing modules on formal  homogeneous spaces.}
\label{sec:dualizing}
 In this section, we assume that  the
$\K$-module $\q=\g/\h$ is free of finite rank $(p,q)$.

 By (\ref{eq:isoq}), the algebra $F(G/H)$ is
formally smooth of dimension $(p,q)$ and the dualizing module
$K(F(G/H))$ is defined. By definition, it is free of rank $1$ as an 
$F(G/H)$-module. For $a\in \g$, $-\zeta_{a}^\lambda$ is a
derivation of $F(G/H)$. It acts in  $K(F(G/H))$, so $K(F(G/H))$ is
in a natural way a $\g-F(G/H)$-module (recall definition
\ref{def:Ind}).  The evaluation
\begin{equation}\label{eq:iso}
\bar \delta: K(F(G/H))
         \to K(F(G/H))\otimes_{F(G/H)} \K
\end{equation}
is a morphism of $\h$-modules. It follows that the
$\K$-isomorphism $K(F(G/H))\otimes_{F(G/H)} \K
\simeq\Ber((\g/\h)^*)$ given in (\ref{eq:canonique4})  is an
isomorphism of $\h$-modules. It follows from proposition
\ref{prop:blattner} that the map $\tilde \delta$ defined in
(\ref{eq:coinducemap}) is an isomorphism
\begin{equation}\label{eq:isob}
\tilde \delta: K(F(G/H)) \to F(G/H, \Ber ((\g/\h)^*)).
\end{equation}
of $\g-F(G/H)$-modules.

For  $a\in \h$,  $\ad a $  induces an endomorphism  of $\g/\h$. We
denote by $\str_{\g/\h}(\ad a )$  its supertrace.  The action of
$a$ in $\Ber \g/\h$ is the multiplication by $\str_{\g/\h}(\ad a
)$, and the action of $a$ in $\Ber ((\g/\h)^*)$ is the
multiplication by $-\str_{\g/\h}(\ad a )$.

\begin{proposition}\label{prop:invariant}
The dualizing module $K(F(G/H))$ has a $\g$-invariant basis as an 
$F(G/H)$-module if and only if the  unimodularity condition
(\ref{eq:unimb}) holds. In this case, the $\K$-submodule  of
$\g$-invariant elements of $K(F(G/H))$ is free of rank one, and a
basis of the $\K$-module of $\g$-invariant is a basis of the
$F(G/H)$-module  $K(F(G/H))$.
\end{proposition}
\begin{proof}
Let $D$ be an invariant basis of $K(F(G/H))$. The element $\bar
\delta (D)$ is a basis of $\Ber (\g/\h)^*$ invariant by $\h$. So
(\ref{eq:unimb}) holds.

Conversely, suppose that (\ref{eq:unimb}) holds. Then $K(F(G/H))$
is isomorphic to $F(G/H)$ --see (\ref{eq:isob})--, and the
constant function $1\in F(G/H)$ is an invariant basis of $F(G/H)$.
\end{proof}

\begin{remark}\label{rem:special}
Similarly, the  $\g-F(G/H)$-module $\CT(F(G/H))$ of derivations of
$F(G/H)$ is isomorphic to $F(G/H, \g/\h)$ and the module
$\Omega_{ev}^1(F(G/H))$ to $F(G/H, (\g/\h)^*)$. However, even
under assumption (\ref{eq:unimb}), there is usually no basis of
$F(G/H, \g/\h)$ consisting  of $\g$-invariant elements. So the
simple proof of proposition \ref{prop:invariant}, given in
proposition \ref{prop:haar} for the special case $\h=0$,  does not
generalize.
\end{remark}

\medskip
To get invariant elements without  unimodularity hypotheses, we
consider twisted dualizing modules. Let $V$ be an $\h$-module.  We
define
\begin{equation}\label{eq:jacobian6}
K(F(G/H),V):=K(F(G/H))\otimes_{F(G/H)} F(G/H,V),
\end{equation}
and we  call it the \emph{dualizing module twisted by $V$}. As in
(\ref{eq:isob}), there is a canonical isomorphism of (recall
(\ref{eq:isoisob}) and (\ref{eq:isob})) $\g - F(G/H)$-modules
\begin{equation}\label{eq:jacobian7}
K(F(G/H),V)\simeq  F(G/H,\Ber ((\g/\h)^*)\otimes_\K V) .
\end{equation}
We apply this construction to the $\h$-module $\Ber \g/\h$. Using
the isomorphism (\ref{eq:canonical8}), we obtain a canonical
isomorphism of $\g - F(G/H)$-modules
\begin{equation}\label{eq:jacobian8}
K(F(G/H),\Ber\g/\h)\simeq   F(G/H).
\end{equation}
We  denote by $\mathbf D_{G/H}$ the element of
$K(F(G/H),\Ber\g/\h)$ which corresponds to $1\in F(G/H)$ by this
isomorphism.

\begin{proposition}\label{prop:canonical}
\label{prop:canonicalBerezinian} i) The element $\mathbf D_{G/H}
\in K(F(G/H),\Ber \g/\h)$ is $\g$-invariant and $\bar \delta
(\mathbf D_{G/H})\in \Ber (\g/\h)^*\otimes_\K \Ber \g/\h$ is equal
to the canonical element $\mathbf d$ defined in
(\ref{eq:canonical9}).

ii) The $\K$-module of $\g$-invariant elements of
$K(F(G/H),\Ber\g/\h)$ is free of rank one with basis $\mathbf
D_{G/H}$.

iii) $\mathbf D_{G/H}$ is a basis of the $F(G/H)$-module
$K(F(G/H),\Ber \g/\h)$.
\end{proposition}

\subsection{The Jacobian of the exponential map on  a formal homogeneous space.}
 \label{sec:jacobians}
We keep the assumptions  of subsection \ref{sec:dualizing}. The
isomorphism $ \widetilde \beta^*:  F(G/H) \to  F(\q) $ of
augmented algebras --see (\ref{eq:isoq})--- induces like in
(\ref{eq:inverse}) a  $\K$-isomorphism $\widetilde \beta^*:
K(F(G/H))\to K(F(\q))$. Since the natural action of $a \in \g$ in
$F(G/H)$ is given by the derivation $-\zeta_{a}$ (see subsection
\ref{sec:32}), the natural action in $K(G/H)$ is given by
$-\CL(\zeta_{a})$. By (\ref{eq:isoq3}) and (\ref{eq:inverse6}),
the natural action of $a$ in $ K(F(\q))$ obtained by transporting
$-\CL(\zeta_{a})$ by $\widetilde \beta^*$ is equal to
$-\CL(\zeta_{\alpha_a}).$

We recall that $\chi$ is a representation of $\h$ in $V$. We
define
\begin{equation}\label{eq:jacobian30}
K(F(\q),V): =K(F(\q))\otimes_{F(\q)} F(\q,V).
\end{equation}
Because of   (\ref{eq:inverse4}), the natural map $\widetilde
\beta^* \otimes \bar \beta^*_\chi: K(F(G/H))\otimes_\K F(G/H,V)
\to K(F(\q),V)$ induces an isomorphism $\bar
\beta^*_\chi:K(F(G/H),V) \to K(F(\q),V)$. Recall that the action
of $a$ in $F(\q,V)$ is denoted by $-\bar \zeta_{\chi,a}$. It
follows that the action of $a$ in $K(F(\q),V)$ is given by
$-\CL(\zeta_{\alpha_a})\otimes 1 - 1 \otimes \bar \zeta_{\chi,a}$.
From (\ref{eq:coderxd}) we get the following formula for the
action of $a\in\g$ in $K(F(\q),V)$:
\begin{equation}
\label{eq:module}
  a(d \otimes \phi)
= -\CL(\zeta_{\alpha_a})(d)\otimes \phi
 \pm
\left(  -d \otimes \zeta_{\alpha_a}(\phi) +  d
\otimes\chi(\theta_\alpha)(\phi) \right) \mbox{ for }  d\in
K(F(\q)), \phi \in F(\q,V)
 .
\end{equation}

We choose  a basis $d_{\q}$ of $K(F(\q))$ which is invariant by
translations (see lemma \ref{lem:haar}). An element of
$K(F(\q),V)$ can be written in an unique manner as $ d_{\q} \phi
$, with $\phi \in F(\q,V)$  ---we leave out the symbol $\otimes$.
Formula (\ref{eq:module}) specializes (see (\ref{eq:divergence})
and (\ref{eq:divergence2}))  to
\begin{equation}\label{eq:module2}
 a( d_{\q}\phi)=
  \pm d_{\q}\left( -\diverg(\zeta_{\alpha_a}) \phi
 -  \zeta_{\alpha_a}(\phi)+ \chi(\theta_a)(\phi)\right)
 \mbox{ for } a\in \g, \phi \in F(\q,V).
\end{equation}

\medskip
We apply these considerations to the $\h$-module $V=\Ber \g/\h$.
We recall (see subsection \ref{sec:dualizing})  that the action of
$\h$ is the multiplication by the linear form $\chi=\str_{\g/\h}$.
We denote by
\begin{equation}\label{eq:canonicald}
\mathbf d_{\q}\in  K(F(\q), \Ber \g/\h)
\end{equation}
the canonical element, that is, recalling the definition
(\ref{eq:canonical9}) the element such that $\delta(\mathbf
d_{\q})=\mathbf d \in \Ber (\g/\h)^*\otimes  \Ber \g/\h$.

An element of $K(F(\q), \Ber \g/\h)$ can be written in a unique
manner as $ \phi \mathbf d_{\q} $, with $\phi \in F(\q)$. From
equation (\ref{eq:module2}) we get
\begin{equation}\label{eq:module22}
 a( \phi \mathbf d_{\q})=
 \left(-
   \diverg(\zeta_{\alpha_a})\, \phi
  -\zeta_{\alpha_a}  (\phi)
  +    \str_{\g/\h}\circ\theta_a\, \phi
 \right) \mathbf d_{\q}
 \mbox{ for } a\in \g, \phi \in F(\q).
\end{equation}

\medskip
We defined  $\mathbf D_{G/H} \in K(F(G/H),\Ber \g/\h)$ in
proposition \ref{prop:canonicalBerezinian}. There exists a
function $J\in F(\q)$ such that
\begin{equation}\label{eq:jacobian5}
\bar \beta^*_{\str_{\g/\h}}(\mathbf D_{G/H})= J \mathbf d_{\q}.
\end{equation}
We call $J$ the \emph{Jacobian of the exponential map $\q \to
G/H$}. A natural problem is to compute $J$. We already did it in
theorem \ref{theo:jac} when $\h=0$. In section
\ref{sec:Berezinians10} we compute $J$ in the  case of symmetric
pairs. We shall use the following characterization of $J$, which
follows immediately from (\ref{eq:module22}).
\begin{lemma}\label{lem:calcul}
The function $J\in F(\q)$ is the unique element of $F(\q)_\0$
which satisfies  $\delta(J)=1$ and such that, for all   $a\in \g$,
\begin{equation}\label{eq:module3}
\frac{  \zeta_{\alpha_a}  (J) }{J} + \diverg(\zeta_{\alpha_a})
-\str_{\g/\h}\circ\theta_a=0.
\end{equation}
\end{lemma}

\begin{remark}
\label{rem:unimod} By proposition \ref{prop:invariant}, if the
unimodularity  condition (\ref{eq:unimb}) holds, there exists a
$\g$-invariant basis $D_{G/H}$ of $K(F(G/H))$. If we ask also for
the equality  $\bar \delta(D_{G/H})=\bar \delta (d_{\q})$   in
$\Ber(\q^*)$, it is unique, and we obtain, with the same function
$J$ than in (\ref{eq:jacobian5}),
\begin{equation}\label{eq:jacobian55}
\widetilde \beta^*(D_{G/H})= J d_{\q}.
\end{equation}
\end{remark}

\subsection{Homogeneous spaces of purely odd dimension.} \label{sec:odd2}
We keep the notations of subsection \ref{sec:dualizing}, and in
addition we assume that $\dim \g/\h=(0,q)$.
We apply the results
of subsection \ref{sec:odd} to the algebra $F(G/H)$.

Since it is free of finite rank, we have
$F(G/H)_*=F(G/H)^*=U(\g)/U(\g)\h$, and Berezin integral provides a
canonical isomorphism $\psi_{F(G/H)}:K(F(G/H)) \to U(\g)/U(\g)\h $
of $\g-F(G/H)$-modules (see subsection \ref{sec:odd}). Let $V$ be
an $\h$-module. Tensoring with $F(G/H,V)$, we obtain a canonical
isomorphism $K(F(G/H),V) \to U(\g)/U(\g)\h \otimes_{F(G/H)}
F(G/H,V)$. Composing with the isomorphism (\ref{eq:isoisoa}), we
obtain a canonical isomorphism of $\g-F(G/H)$-modules
\begin{equation}\label{eq:isoisoca}
\kappa_{\chi}:  K(F(G/H),V)  \simeq U(\g)\otimes_\h V,
\end{equation}
and, using (\ref{eq:jacobian7}), a canonical isomorphism of
$\g-F(G/H)$-modules
\begin{equation}\label{eq:isoisoc}
F(G/H,\Ber(\g/\h)^*\otimes V) \simeq U(\g)\otimes_\h V.
\end{equation}

\begin{remark}\label{rem:classique}
The isomorphism (\ref{eq:isoisoc}) between an induced and a
coinduced representation of $\g$ is well known, see
\citelist{\cite{Chemla1993} \cite{Chemla1994} \cite{Chemla1994b}
\cite{BellFarnsteiner1993}}. It is used in Gorelik's paper
\cite{Gorelik2001}. The properties we need in this article, that
is the relation with dualizing modules and the fact that it is an
isomorphism of $F(G/H)$-modules, are explicit in
\citelist{\cite{Chemla1993} \cite{Chemla1994} \cite{Chemla1994b}}.
\end{remark}
In particular, using the canonical isomorphism
(\ref{eq:canonical8}), the $\g-F(G/H)$-module $U(\g)\otimes _\h
\Ber \g/\h$ is canonically isomorphic to $F(G/H)$. We denote by
$\mathbf T \in U(\g)\otimes _\h \Ber \g/\h$ the element
corresponding to $1\in F(G/H)$. Remember (see proposition
\ref{prop:canonicalBerezinian}) the definition of the canonical
element $\mathbf D_{G/H} \in K(F(G/H),\Ber \g/\h)$. We have
\begin{equation}\label{eq:kappaD}
\kappa_{\str_{\g/\h}}(\mathbf D_{G/H})=\mathbf T.
\end{equation}
We have proved the following proposition, analogous to proposition
\ref{prop:canonicalBerezinian}.

\begin{proposition}\label{prop:classique}
a)  The $\K$-module of $\g$-invariants of  $U(\g)\otimes _\h
\Ber\g/\h$ is free of rank~$1$ with basis $\mathbf T$. Moreover,
$\mathbf T$ is a basis of the $F(G/H)$-module $U(\g)\otimes _\h
\Ber \g/\h$, and the map $F(G/H)\ni \phi \to \mathbf T\phi \in
U(\g)\otimes _\h \Ber \g/\h $ is an isomorphism of
$\g-F(G/H)$-modules.

b) Assume  condition (\ref{eq:unimb}). Then the $\K$-module of
$\g$-invariants of  $U(\g)/U(\g)\h $  is free of rank~$1$.    A
basis $T$ of this module  is also  a basis of the  $F(G/H)$-module
$U(\g)/U(\g)\h $. The map $\phi \to
  T\phi$ from $F(G/H)$ to  $U(\g)/U(\g)\h $ is an
isomorphism of $\g-F(G/H)$-modules.
\end{proposition}

\bigskip
Similarly,  we have the isomorphism of $\CT(F(\q))-F(\q)$-modules
$\psi_{F(\q)}: K(F(\q)) \to S(\q)$, and lemma \ref{lem:dualite2}
gives
\begin{equation}\label{eq:dualite3}
\psi_{F(G/H)}=\widetilde \beta\circ \psi_{F(\q)}\circ\widetilde
\beta^*.
\end{equation}
Tensoring with $F(\q,V)$, we obtain a canonical isomorphism
\begin{equation}\label{eq:dualite9}
\bar \kappa_{\chi} : K(F(\q),V) \to
S(\q)\otimes_{F(\q)}F(\q,V)\simeq S(\q)\otimes V.
\end{equation}
Formula (\ref{eq:dualite3}) extends to
\begin{equation}\label{eq:dualite12}
\kappa_{\chi}=\bar \beta_{\chi}   \circ \bar \kappa_{\chi}
 \circ\bar \beta^*_\chi.
\end{equation}

We identify $\q^*$ and $(\g/\h)^*$. Since the rank of $\g/\h$ is
$(0,q)$, we have $\Ber \g/\h=S^q(\q^*)$. There is a canonical
element $\mathbf t \in S(\q)\otimes\Ber \g/\h $: in term of a
totally ordered basis $(e_i)_{1\leq 1\leq q}$ of $\q$ and of its
dual basis $(x^i)_{1\leq i\leq q}$ of $\q^*$, we have
\begin{equation}\label{eq:TTTTT}
\mathbf t = (-1)^q\,  e_1\cdots e_q\otimes x^q\cdots x^1.
\end{equation}
Recall the definition (\ref{eq:canonicald}) of $\mathbf d_\q\in
K(F(\q),\Ber \g/\h)$. The analogue of (\ref{eq:kappaD}) is the
formula
\begin{equation}\label{eq:barkappaD}
\bar \kappa_{\str_{\g/\h}}(\mathbf d_{\q})=\mathbf t.
\end{equation}

\medskip
Recall  the Jacobian  $J\in F(\q)$ of the exponential map $\q \to
G/H$ (formula (\ref{eq:jacobian5})). {\it The following theorem
explains the role of this Jacobian
 in the study of Gorelik's elements}.
\begin{theorem}\label{thm:correspondance}
\begin{enumerate}
\item[a)] We have $\bar \beta_{\str}(J \mathbf t)=\mathbf T$.
\item[b)] Suppose  that condition  (\ref{eq:unimb}) holds. Then
$\widetilde \beta(J \, e_1\cdots e_q) \in  U(\g)/U(\g)\h $ is a
basis of the $\K$-module of  $\g$-invariants of  $U(\g)/U(\g)\h $,
and a basis of  $U(\g)/U(\g)\h $ as a $F(G/H)$-module.
\end{enumerate}
\end{theorem}
\begin{proof} a) From (\ref{eq:kappaD}) and  (\ref{eq:jacobian5})
we get
$$
\bar \beta^*_{\str_{\g/\h}}(\kappa_{\str_{\g/\h}}^{-1}(\mathbf
T))= J \mathbf d_{\q}.$$ Using the fact that $\bar
\kappa_{\str_{\g/\h}}$ is a $F(\q)$-module isomorphism, we obtain
$$
\bar \kappa_{\str_{\g/\h}}(\bar
\beta^*_{\str_{\g/\h}}(\kappa_{\str_{\g/\h}}^{-1}(\mathbf T)))= J
\bar \kappa_{\str_{\g/\h}}(\mathbf d_{\q}).
$$
We get our formula using (\ref{eq:barkappaD}) and
(\ref{eq:dualite12}).

 b) Follows from part a) and proposition \ref{prop:classique}.
\end{proof}

\section{Symmetric pairs.}
\label{sec:symmetric} In this section, we consider a Lie symmetric
pair  $(\g,\sigma)$ with corresponding decomposition $\g=\h\oplus
\q$.

\subsection{Representations by coderivations.} \label{sec:coderivations}
Let $a\in \g$. Recall
the coderivation $\widetilde l_a=c_{\alpha_a}$ of $S(\q)$ (formula
(\ref{eq:isoq3})). In this subsection, we give a formula for the
vector field $\alpha_a\in F(\q,\q)$.  In fact, we do more.
Following the method of \cite{Petracci2003}, we describe \emph{all
the representations by coderivations of $\g$ in $S(\q)$ which are
universal}. Informally, this means that they are given by formulas
which are the same for all symmetric pairs $(\g,\sigma)$.
Formally, it can be express in functorial terms (see
\cite{Petracci2003}). Reasoning as in \cite{Petracci2003}, these
representations are of a very special form, and we are led to the
following considerations.

\medskip
We denote by $y:=x_\q\in F(\q,\q)$ the generic point of $\q$. We
consider $y$ as the element of the Lie $F(\q)$-superalgebra
$F(\q,\g)$ which is $0$ on $S^n(\q)$ for $n\neq 1$, and which is
the natural injection $\q\to \g$ on $\q$. The decomposition
$F(\q,\g)=F(\q,\h)\oplus F(\q,\q)$ is a symmetric pair in the
category of Lie $F(\q)$-superalgebras. For a formal power series
$p(t) \in\K_0[[t]]$ the element $p(\ad y)\in
\End_{F(\q)}(F(\q,\g))$ is defined by   locally finite sums  like
in (\ref{eq:sum}). In particular, for  $a \in \g$, the element
$p(\ad y)(a)\in F(\q,\g)$ is defined  by the following formula,
similar to (\ref{eq:radx}). Let $b_1$, $b_2$, \dots, $b_n$ in
$\q$. We have
\begin{equation}\label{eq:radx2}
        p(\ad y)(a) : b_1\cdots b_n
\mapsto p_n
        \sum_{s\in S_n}
        \pm \ad  {b_{s(1)}}\cdots \ad  {b_{s(n)}}  (a).
\end{equation}

We consider two  series $d(t) , p(t) \in\K_0[[t]]$ such that
$p(t)=p(-t)$ and $d(t)=-d(-t)$. We write $ p(t)=p_0+p_2t^2+\cdots
\in \K_\0[[t^2]]$   and $d(t)=d_1t+d_3t^3+\cdots \in
t\K_\0[[t^2]].$ Formula (\ref{eq:radx2}) implies that for $a\in
\q$ we have $p(\ad y)(a) \in F(\q,\q)$, and for $a\in \h$ we have
$d(\ad y)(a)\in F(\q,\q)$. In particular, the following definition
is meaningful.
\begin{definition}\label{def:ca}
For  $a\in \q$, we denote by $C^a$ (or $C^a_{p,d,\g,\sigma}$ if we
need to specify) the coderivation $c_{p(\ad y)(a)}$ of $F(\q)$.
For $a\in \h$, we denote by  $C^a$ the coderivation $c_{d(\ad
y)(a)}$ of $F(\q)$. For $a\in \g$, we define $C^a$ by additivity.
\end{definition}

The bracket $[C^a,C^b]$ has been computed in \cite{Petracci2003}.
Reasoning as in  \cite{Petracci2003}, we obtain the following
theorem.

\begin{theorem}\label{theorem:emanuela}
The map $a \to C^a_{p,d,\g,\sigma}$ is a representation of $\g$ in
$S(\q)$ for any symmetric pair $(\g,\sigma)$ if and only if the
following functional equations hold in $\K[[t,u]]$:
\begin{eqnarray}
\label{pippo1}
   d(u)\frac{d(t+u) - d(t)}{u}
 + d(t)\frac{d(t+u) - d(u)}{t}
 = -d(t+u),
\\
\label{pippo2}
   p(u)\frac{p(t+u) - p(t)}{u}
 + p(t)\frac{p(t+u) - p(u)}{t}
 = -d(t+u),
\\
\label{pippo3}
   p(u)\frac{d(t+u) - d(t)}{u}
 + d(t)\frac{p(t+u) - p(u)}{t}
 = -p(t+u).
\end{eqnarray}
\end{theorem}

Even if we know it in advance, it is remarkable that these
functional equations have non trivial solutions.

\begin{lemma}\label{lem:zero}
Let $(p,d)$ be a solution of the equations of theorem
\ref{theorem:emanuela}.

If  $d=0$ or $p=0$, then  $(p,d)=(0,0)$.

Assume $d\neq 0$. If the first non zero coefficient  of $d$ is not
a divisor of $0$ in $\K_\0$,  then $d=-t$ and  $p_0\neq 0$.
\end{lemma}
\begin{proof}
The first assertion follows from formulas (\ref{pippo3}) and
(\ref{pippo2}).

Formula (\ref{pippo1}) evaluated at $u=0$ gives $ d(t)
\left(\frac{d(t)}{t}+1 \right)=0 $. If the first non zero
coefficient  of $d$ is not a divisor of $0$ in $\K_\0$,  the
multiplication by $d$ is injective, and we get $ \frac{d(t)}{t}+1
=0 $.

We now assume that $d=-t$. Equations (\ref{pippo1}) and
(\ref{pippo3}) are satisfied, and we are left with the functional
equation
\begin{equation}\label{eq:fonctionnelle}
p(u)\frac{p(t+u) - p(t)}{u} + p(t)\frac{p(t+u) - p(u)}{t} =t+u.
\end{equation}
Evaluating at $u=0$ we obtain
\begin{equation}\label{eq:differentielle}
p_0 p'(t) + p(t)\frac{p(t) - p_0}{t} =t.
\end{equation}
Suppose  $p_0=0$. This gives $p(t)^2=t^2$. Since $p(t)\in
t^2\K[[t^2]]$, this is impossible.
\end{proof}

Lemma \ref{lem:zero} implies that all the non zero solutions of
the functional equations of theorem \ref{theorem:emanuela} with
coefficient in a field $k\subset \K_\0$ are given by  theorem
\ref{theorem:emanuela2} below. We recall that the Bernoulli
numbers $b_n\in \Q$ are defined by the generating series
$\frac{t}{e^t-1} =  \sum_{n=0}^\infty b_{n} \frac{ t^{ n}}{  n
!}$. For $c \in \K_\0^\times$, we define $p_c(t) \in \K_\0[[t^2]]$
by
\begin{equation}\label{eq:differentielle4}
p_c(t)= t\coth\left(\frac{t}{c}\right)=c+ \sum_{n=1}^\infty b_{2n}
\frac{ 2^{2n}t^{2n}}{c^{2n-1}(2n)!}= c + \frac{1}{3c}t^2 -
\frac{1}{45c^3}t^4 +\cdots.
\end{equation}
\begin{theorem}\label{theorem:emanuela2}
Let $c \in \K_\0^\times$. There exists an unique element $p(t)\in
\K_\0[[t^2]]$ such that $p(0)=c$ and such that the pair
$(p(t),-t)$ is a solution of the functional equations of theorem
\ref{theorem:emanuela}.  We have $p(t)=p_c(t)$.
\end{theorem}
\begin{proof}
We repeat the arguments of \cite{Petracci2003}. The pair
$(p(t),-t)$ is a solution if and only if the pair
$(p(\frac{t}{c})c , -t)$ is a solution. So we assume  that
$p_0=1$. In this case, according to formula
(\ref{eq:differentielle}), $p$ will be a solution of the
(singular) Cauchy problem
\begin{equation}\label{eq:differentielle2}
t p'(t) + p(t)(p(t) - 1) =t^2 \mbox{  and   } p(0)=1.
\end{equation}
Writing (we  forget for a moment that $p$ is supposed to be even)
$p(t)=1+p_1 t + p_2 t^2+ \dots$, we get inductively $p_1=0$,
$p_2=\frac{1}{3}$, and $p_n$ in function of the $p_k$ with $k<n$.
This proves that there is at most one solution of
(\ref{eq:differentielle2}).

We leave to the reader to check that the function $t\coth(t)$ is
indeed a solution of the differential equation
(\ref{eq:differentielle2}), and --- this is much stronger--- of
the functional equation (\ref{eq:fonctionnelle}).
\end{proof}

\begin{remark}\label{rem:carN}
As in  \cite{Petracci2003}, theorem \ref{theorem:emanuela} is
meaningful and true without assuming that $\K_\0$ contains $\Q$,
but only assuming that $2$ is invertible in $\K$. However, the
functional equations have   no solutions if $\K_\0$ does not
contains $\Q$.
\end{remark}

Let $c\in \K_\0^\times$. We denote by $a \to \alpha_c^a$ the
linear map from $\g$ to the space of $F(\q,\q)$ of vector fields
on $\q$ such that
\begin{equation}\label{eq:vectorfield}
 \alpha_c^a= [a,y] \mbox{ for } a\in \h, \mbox{ and }
  \alpha_c^a= p_c(\ad y)(a) \mbox{ for } a\in \q.
\end{equation}
We denote by
\begin{equation}\label{eq:vectorfield4}
C_c^a:=c_{\alpha_c^a}
\end{equation}
the corresponding coderivation of $S(\q)$. It follows from theorem
\ref{theorem:emanuela} and theorem \ref{theorem:emanuela2} that
$C_c: a\to C^a_c$ is a representation of $\g$ in $S(\q)$.

\begin{remark}\label{rem:equivalent}
We denote by $I_c$ the automorphism of the Hopf algebra $S(\q)$
which is the multiplication by $c^n$ in $S^n(\q)$. It intertwines
the representations  $ C_1 $ and $C_c$ in $S(\q)$. Thus, if $\K$
is a field of characteristic $0$,   there are, up to equivalence,
exactly two universal representations by coderivations of $\g$ in
$S(\q)$: the trivial one and $  C_1 $.
\end{remark}

\begin{remark}\label{rem:deriv}
In term of representations, the equality $d=-t$ means the
following. Let $a\in \h$. Then $C_c^a$ (which does not depend on
$c$) is a derivation of $S(\q)$ (of course, by definition, it is
also a coderivation). It is determined by the fact that
\begin{equation}\label{eq:Caa}
C_c^a(b)= [a,b] \mbox{  for   } b\in \q.
\end{equation}
Thus, $C^a_c$ is the natural adjoint action in of $a\in \h$ in
$S(\q)$.
\end{remark}

For $a\in \q$ and $w \in  S ^n(\q)$, formula
(\ref{eq:differentielle4}) gives an explicit expression  for
$C_c^a(w)\in S(\q)$  in terms of the Bernoulli numbers and of the
coproduct  of $w$. We give a few terms. For $b,b_1,b_2 \in \q$, we
have
\begin{equation}\label{eq:Caab}
      C_c^a(1)= c\, a,
\quad C_c^a(b)= c\, ab, \quad C_c^a(b_1b_2) = c\, ab_1b_2 +
\frac{1}{3c}\Big( [[a,b_1],b_2]\pm[[a,b_2],b_1] \Big) ;
\end{equation}
and  for $n>2$ and  $b_1,\dots,b_n \in \q$, we have
\begin{equation}\label{eq:Caab7}
  C_c^a(b_1 \cdots b_n)
= c\, ab_1 \cdots b_n +\frac{1}{3c}
  \sum_{1\leq i< j \leq n}
  \Big(\pm [[a,b_i],b_j]\pm[[a,b_j],b_i]
  \Big)
  b_1 \cdots \hat b_i \cdots \hat b_j \cdots b_n
+ \dots
\end{equation}
where the hat symbol means omission, and the non written terms are
of degree $< n-2$.

\subsubsection{Formulas for and $\widetilde l_a$ and $\alpha_a$}
The map $a \to C_c^a$ extends to an homomorphism of algebras $u\to
C_c^u$ from  $U(\g)$ to $\End(S(\q))$. For $u\in U(\g)$ we write
$\tilde \tau(u)=C_1^u(1)$. This defines a map $\tilde \tau : U(\g)
\to S(\q)$ which intertwines the left action of $\g$ in $U(\g)$
and the representation $C_1$. Since $C_1^a(1)=0$ for $a\in \h$,
this gives an intertwining  map
\begin{equation}\label{eq:tau}
\tau :U(\g)/U(\g)\h \to S(\q).
\end{equation}
The map $\tau$ can be explicitly written in terms of Bernoulli
numbers. For example, for $b,b_1,b_2,b_3$ in $\q$ we have
\begin{equation}\label{eq:tau1}
\tau(1)=1,\quad \tau(j(b))=b, \quad \tau(j(b_1)j(b_2))=b_1b_2,
\end{equation}
\begin{equation}\label{eq:tau2}
 \tau(j(b_1)j(b_2)j(b_3))
  =b_1b_2b_3
  + \frac{1}{3}\Big([[b_1,b_2] ,b_3]\pm [[b_1,b_3], b_2]\Big).
\end{equation}

For the following theorem, recall the notations $\alpha_a$ and
$\widetilde l_a$ of proposition \ref{prop:formulacoder} and the
notation $\alpha_1^a$ of (\ref{eq:vectorfield}).
\begin{theorem}\label{theorem:emanuela3}
i) The map $\tau $ is the inverse of the map $\widetilde \beta
:S(\q) \to U(\g)/U(\g)\h $.

ii) For $a\in \g$, we have $\widetilde  l_a=C^a_1$

iii) For $a\in \g$, we have $\alpha_a=\alpha^a_1$.
\end{theorem}
\begin{proof}
The proof of part i) is similar to the proof of
Poincar{\'e}-Birkhoff-Witt theorem in \cite{Petracci2003}. Here is a
brief sketch. Let $b \in \q_\0$ and $n \in \N$.  By induction on
$n$, we obtain that $\tau(j(b)^n)=b^n$. Applying this to the
generic point $y$ of $\q$ in the Lie superalgebra $F(\q,\g)$
proves part a). Part  ii) follows from part i), and part iii) is a
rephrasing of ii).
\end{proof}

\begin{remark} \label{rem:koszul}
When $\g=\g_\0\oplus \g_\1$ is a Lie superalgebra over a field of
characteristic zero, and $\h=\g_\0$, the formula $\widetilde
l_a=C^a_1$ is due to Koszul \cite{Koszul1982}. Presumably, his
proof extends to the general case treated in theorem
\ref{theorem:emanuela3}. However, we think that our method of
proof is interesting: it says that the properties of Bernoulli
numbers which make the theorem true are encoded by the fact that
$t\coth(t)$ satisfies the functional equation
(\ref{eq:fonctionnelle}).
\end{remark}

\subsection{Induced and coinduced representations.}\label{sec:universal2}
We keep the   notations of subsection \ref{sec:coderivations}. We
consider in addition a representation $\chi: \h \to \End(V)$ of
$\h$ in a $\K$-module $V$.

Let $a\in \g$. We  give an explicit  formula for $\bar l_{\chi,a}
$ defined in (\ref{eq:coderxc}). We follow a method similar to the
one of subsection \ref{sec:coderivations}, where we consider the
case of the trivial representation of $\h$ in $\K$.

Let   $(h,q)$ be a pair of formal series $h=h_0+h_2 t^2+\cdots \in
\K_\0[[t^2]]$   and $q(t)=q_1t+q_3 t^3 +\cdots \in t\K_\0[[t^2]]$.
We define a linear map   $\g\ni a  \to \theta^a\in F(\q,\h)$  by
the  formulas
\begin{equation}\label{eq:theta1}
\theta^a= h(\ad y)(a) \mbox{ for } a\in \h, \mbox{ and }\theta^a=
q(\ad y)(a) \mbox{ for } a\in \q.
\end{equation}

Let $c \in \K_\0^\times$, and recall   the representation $C_c$ of
$\g$ in $S(\q)$ (see (\ref{eq:vectorfield4})). Following the model
(\ref{eq:coderxc}), we define an endomorphism $\Theta_{\chi}^a \in
\End(S(\q)\otimes V)$ (or $\Theta_{c,\chi,h,q,\g,\sigma}^a$ if we
need to specify) by the following formula:
\begin{equation}\label{eq:theta1a}
  \Theta_{\chi}^a(w\otimes v)
= C^a_c(w)\otimes v + (\id_{S(\q)}\otimes\chi\circ\theta^a) \circ
\Delta_{S(\q)}(w)(1\otimes v) .
\end{equation}
For example, for $a\in \h$, we have:
\begin{equation}\label{eq:theta2}
\Theta_{\chi}^a(1\otimes v) =1\otimes \chi( h_0 a) (v).
\end{equation}
We denote by $\Theta_{\chi}$  the map $a \to \Theta_{\chi}^a$, and
we study when it is a representation.

\begin{lemma} \label{lem:rep}
i) The map $\Theta_{\chi} $ is a representation of $\g$  if and
only if
$$
\chi\circ(  \theta^a\circ C^b_c - \pm
  \theta^b\circ C_c^a
+ [\theta^a,\theta^b] - \theta^{[a,b]} ) = 0
 \ \ \forall a,b\in\g
.
$$

ii) The map $\Theta_{\chi} $ is a representation of $\g$ for all
representations $\chi$ of $\h$  if and only if
$$
 \theta^a\circ C^b_c - \pm
  \theta^b\circ C^a_c
+ [\theta^a,\theta^b] - \theta^{[a,b]} = 0
 \ \ \forall a,b\in\g
.
$$
\end{lemma}
\begin{proof}
We write $\Theta=\Theta_\chi$, $C=C_c$, and
$\Delta=\Delta_{S(\q)}$ for simplicity.

$i)$ Let $w\otimes v\in S(\q )\otimes V$ and $a,b\in\g $. As $C^a$
is a coderivation we get
\begin{eqnarray*}
  &&\Theta^a\circ \Theta^b
    (w\otimes v)=
\\
 &=& C^a\circ C^b(w)\otimes v
   + \left( C^a\otimes \chi \circ\theta^b
     \pm
      C^b\otimes \chi \circ \theta^a
     \right)
\circ\Delta(w)(v)+
\\
&&  + (\id\otimes \chi\circ  \theta^a \circ C^b
      )
\circ \Delta(w)(v)
    + (\id\otimes \chi \circ \theta^a \otimes
            \chi \circ\theta^b
      )
\circ (\Delta\otimes \id) \circ \Delta(w)(v) .
\end{eqnarray*}
Since $\chi$ is a representation we have $      \left(
\chi\circ\theta^a
        \otimes \chi\circ\theta^b
         - \pm
           \chi\circ\theta^b
         \otimes  \chi\circ\theta^a
      \right)
\circ \Delta
    = \chi\circ[\theta^a,\theta^b]$.
Since  $\Delta$ is associative  and $C$   a representation  we
have
\begin{eqnarray*}
  &&[\Theta^a, \Theta^b]
    (w\otimes v)=
\\
&=& C^{[a,b]}(w)\otimes v + \left( \id\otimes \chi\circ
         (\theta^a\circ C^b
           -\pm
        \theta^b\circ C^a
         )
  \right)
\circ \Delta(w)(v)+
\\
&&+ \left( \id
   \otimes \chi\circ[\theta^a,\theta^b]
\right) \circ \Delta(w)(v) .
\end{eqnarray*}
This formula shows that  $\Theta$ is a representation if and only
if
\begin{equation} \label{eq:newbis}
\left(
 \id
\otimes
  \chi\circ\left(
          \theta^a\circ C^b
        - \pm
          \theta^b\circ C^a
        + [\theta^a,\theta^b]
        - \theta^{[a,b]}
  \right)
\right) \circ \Delta = 0 \ \  \ \ \forall a,b\in\g .
\end{equation}
To show that the direct part assertion of the lemma is a necessary
condition, it is enough to compose  $\delta\otimes \id_{U(\g)}$
with the left hand side of (\ref{eq:newbis}).

$ii)$ To show the direct part is   sufficient to consider a
faithful representation  $\chi$ of $\h $, for instance the left
regular representation in $U(\h)$.
\end{proof}

\begin{theorem} \label{tm24}
The   map $\Theta_{c,\chi,h,q,\g,\sigma}$ is a  representation for
any symmetric pair $(\g,\sigma)$ and any representation $\chi$ of
$\h$, if and only if the functional equations
\begin{equation}\label{eq7}
  - h(t+u) + h(t) + h(u) - h(t)h(u) = 0,
\end{equation}
\begin{equation}\label{eq8}
\ \frac{q(t+u) - q(t)}{u}
  p_c(u)
+ \frac{q(t+u) - q(u)}{t}
  p_c(t)
- q(t)q(u)+ h(t+u) = 0,
\end{equation}
\begin{equation}\label{eq9}
  q(t) + \frac{h(t+u) - h(u)}{t}
  p_c(t)
- q(t)h(u) = 0,
\end{equation}
are satisfied in $\K_\0[[t,u]]$.
\end{theorem}
\begin{proof}
Let $a$ and $b$ in $\g$. For $X \in \g_\0$, $i,j \in \N$, we set
$$
(t^iu^j:[a,b])_X=[(\ad X)^ia, (\ad X)^j b].
$$
We use the same notation for the Lie algebra $F(\q,\g)$. Thus, for
$a$ and $b$ in $\g$, and the generic point $y\in F( \q, \g)$ we
get an element $(t^iu^j:[a,b])_y \in F(\q,\g)$.  Let $k(t,u) =\sum
k_{ij} t^iu^j\in \K_\0[[t,u]]$. We define
$$
(k(t,u):[a,b])_y:= \sum k_{ij}(t^iu^j:[a,b])_y \in F(\q,\g).
$$
This is possible because this  sum is locally  finite.

Let us denote $h^b(t) = h(t)$ if $b\in\h $, $h^b(t)= q(t)$ if
$b\in\q $, $p_c^b(t) = p_c(t)$ if $b\in\q $, and $p_c^b(t)= -t$ if
$b\in\h $. We have (see \cite{Petracci2003})
\begin{eqnarray*}
&&  \theta^a\circ C^b_c  = \left( - \frac{h^a(t+u)-h^a(t)}{u}
  p^b_c(u)
: [a,b]
  \right)_y=
  \left(
 \frac{h^a(t+u)-h^a(u)}{t}
  p^b_c(t)
: [b,a]
  \right)_y
\\
&&  \left[\theta^a, \theta^b
  \right]
= (h^a(t)h^b(u) : [a,b])_y
\\
&&  \theta^{[a,b]} = ( h^{[a,b]}(t+u): [a,b])_y .
\end{eqnarray*}
Lemma \ref{lem:rep} shows that  $\Theta_\chi $ is a representation
of $\g$ for all representations $\chi$ of $\h$   if and only if
\begin{eqnarray*}
\left(
   h(t)
 - h(t+u)
 + h(u)
 - h(t)
   h(u)
: [a,b] \right)_y = 0 \ \ \ \forall a,b\in\h ,
\\
\Big(
  \frac{q(t+u)- q(t)}{u}
   p_c(u)
 + \frac{q(t+u)-q(u)}{t}
   p_c(t)+
\quad\quad\quad\quad\quad\quad\quad\quad\quad\quad
\\
\quad\quad\quad\quad\quad\quad\quad\quad\quad\quad
 - q(t)
   q(u)
 + h(t+u)
: [a,b] \Big)_y=0 \ \ \forall   a,b\in\q ,
\\
\left(
  \frac{h(t+u)-h(t)}{u}
   p_c(u)
 + q(u)
 - h(t)
   q(u)
: [a,b] \right)_y = 0 \ \ \ \forall   a\in\h,\ b\in\q ,
\\
\left(
   q(t)
 + \frac{h(t+u)-h(u)}{t}
   p_c(t)
 - q(t)
   h(u)
: [a,b] \right)_y = 0 \ \ \ \forall   a\in\q,b\in\h .
\end{eqnarray*}
The last two equations are identical. As in \cite{Petracci2003} we
see that these equations are verified for any Lie symmetric pair
$(\g,\sigma)$ if and only if the functional equations in the
statement of the theorem are satisfied.
\end{proof}

\medskip
We study the  functional equations of Theorem \ref{tm24}.
\begin{lemma}
i) Equation (\ref{eq7}) has an unique solution $h(t)\in
\K_\0[[t^2]]$ such that $h(0)=0$. It is $h=0$.

ii) Equation (\ref{eq7}) has an unique solution $h(t)\in
\K_\0[[t^2]]$ such that $h(0) $ is invertible. It is $h=1$.
\end{lemma}
\begin{proof}
$i)$ If $h(0)=0$ we can divide the equation (\ref{eq7}) by $t$
getting $
  \frac{-h(t+u)+h(u)}{t}
= \frac{h(t)}{t}
  (h(u)-1)
$. The evaluation at $t=0$ gives $ -h'(u)=h'(0)(h(u)-1) . $ As
$h(u)=h(-u)$ we have $h'(0)=0$, so $h'(u)=0$ and  $h(u)$ is
constant. This proves $i)$.

$ii)$ Evaluating at $t=0$ we get $
  0
= h(0)(h(u)-1) . $ If $h(0)$ is invertible, the unique solution is
$h(u)=1$.
\end{proof}

\medskip
\begin{remark}
If $h=0$, equation (\ref{eq9}) gives $q=0$. In this case,
$\Theta_{\chi}$ does not depend on $\chi$. It  is the tensor
product of the representation $C_c$ in $S(\q)$ with the trivial
representation of $\g$ in the module $V$.
\end{remark}

\medskip
We now consider the more interesting case $h=1$. Equations
(\ref{eq7}) and (\ref{eq9}) are satisfied, and  equation
(\ref{eq8}) becomes
\begin{equation}\label{eq:functq}
\frac{q(t+u) - q(t)}{u} p_c(u) + \frac{q(t+u) - q(u)}{t}p_c(t)
  =q(t)q(u)- 1.
\end{equation}

We introduce
\begin{equation}\label{eq:qc}
 q_c(t)
 = -\th\left( \frac{t}{2c}
       \right)
 =-2 \sum_{n=1}^\infty  b_{2n} \frac{
(2^{2n}-1) t^{2n-1}}{c^{2n-1}(2n)!}= -\frac{1}{2c}t +
\frac{1}{24c^3}t^3 +\cdots.
\end{equation}
\begin{remark}\label{rem:th}
The Taylor series of $\th(t)$ is of course well known. However,
since it plays a role in the proof of theorem \ref{thm:main}
below, we remark that  (\ref{eq:qc})  follows from
(\ref{eq:differentielle4}) and the relation
\begin{equation}\label{eq:qccc}
q_c(2t) =\frac{p_c(t)-p_c(2t)}{t},
\end{equation}
which is equivalent to the identity
$\th(t)+\coth(t)-2\coth(2t)=0$.
\end{remark}

\begin{proposition}\label{prop:tanh}
The function $q_c$ is the unique  solution $q\in t\K_\0[[t]]$  of
the functional  equation (\ref{eq:functq}).
\end{proposition}
\begin{proof}
Since $q(0)=0$, the evaluation at $u=0$ of  (\ref{eq:functq})
gives
$$
 q'(t) c + \frac{q(t)}{t} p_c(t)\equiv q'(t) c + q(t) \coth(\frac{t}{c})=-1.
$$
This differential equation has exactly  one solution such that
$q(0)=0$, because the equation recursively determines all the
coefficients of $q$. For instance, for the first one we obtain $
q_1 c + q_1 c=-1$  and  $q_1=-\frac{1}{2c}$, which is compatible
with (\ref{eq:qc}). This solution is of course $q_c$, and one
checks that $q_c$ also satisfy the functional  equation
(\ref{eq:functq}).
\end{proof}

We summarize the main result of this section.  Let $c\in
\K_\0^\times$. We denote by $\g\ni a \to \theta_c^a\in F(\q,\h)$
the linear map  such that
\begin{equation}\label{eq:thetac}
 \theta_c^a= a \mbox{ for } a\in \h, \mbox{ and }
  \theta_c^a=  q_c(\ad y)(a) \mbox{ for } a\in \q.
\end{equation}
We use the notation
$\Theta_{c,\chi}=\Theta_{c,\chi,1,q_c,\g,\sigma}$ (see formula
(\ref{eq:theta1a})). For $a\in \h$, $\Theta_{c,\chi}^a$ does not
depend on $c$. It is the natural endomorphism of $S(\q)\otimes V$
coming from the adjoint action of $a$  in $S(\q)$ and from the
action $\chi(a)$ in $V$. For $a\in \q$, we have
\begin{equation}
\label{eq:theta10} \Theta^a_{c,\chi}  =C^a_c \otimes \id_V +
\left(\id_{S(\q)}\otimes \chi\circ\theta^a_c\right) \circ
\Delta_{S(\q)} .
\end{equation}

\begin{corollary}\label{cor:representation}
i) $\Theta_{c,\chi}$ is a representation of $\g$ in $S(\q)\otimes
V$.

ii) $I_c\otimes \id_V$ intertwines $\Theta_{1,\chi}$ and
$\Theta_{c,\chi}$.
\end{corollary}

\medskip
We use  the notations of the beginning of this section. The
following theorem contains theorem \ref{theorem:emanuela3} as a
particular case. We recall the notations $\alpha_a$, $\alpha_1^a$,
$\theta_a$, $\theta_1^a$ from proposition \ref{prop:formulacoder}
and formulas (\ref{eq:vectorfield}), (\ref{eq:thetac}).
\begin{theorem}
\label{theo:formula} i) For $a\in \g$, we have $\bar
l_{\chi,a}=\Theta_{1,\chi}^a$.

ii) For $a\in \g$, we have $\alpha_a=\alpha_1^a$ and
$\theta_a=\theta_1^a$.
\end{theorem}
\begin{proof}
i) We extend $\Theta_{1,\chi}$ to a map from $U(\g)$ to
$\End(S(\q)\otimes V)$. Since $h=1$ it induces the well defined
map map $\tau : U(\g)\otimes_\h V\ni u\otimes v \to
\Theta_{1,\chi}^u(1\otimes v)\in S(\q)\otimes V$ (recall formula
(\ref{eq:theta2})). By definitions, $\tau$ intertwines the  left
multiplication in $U(\g)\otimes_\h V$ and the representation
$\Theta_{1,\chi}$ in $ S(\q)\otimes V$. As in the proof of theorem
\ref{theorem:emanuela3}, we see that $\tau \circ (\beta\otimes
\id_V)$ is the identity map of $S(\q)\otimes V$.

ii) It is a rephrasing of i).
\end{proof}

\begin{remark} \label{rem:koszul2}
When $\g=\g_\0\oplus \g_\1$ is a Lie superalgebra over a field of
characteristic zero, and $\h=\g_\0$, the formula $\bar
l_{\chi,a}=\Theta^a_{1,\chi}$ is   due to Koszul
\cite{Koszul1982}. For real analytic  symmetric space, the
function $\th(\frac{t}{2})$ occurs in \cite{Rouviere1994}.
\end{remark}

\medskip
We get a similar result for the coinduced representation
$F(G/H,V)$, and the corresponding representation in $F(\q,V)$. Indeed:
from (\ref{eq:coderxd}) and theorem \ref{theo:formula} we obtain
the following   corollary.
\begin{corollary}\label{cor:formula}
For $a\in \g$, we have
\begin{equation}\label{eq:coderxdbis}
 -\bar \zeta_{\chi,a}=-\zeta_{\alpha_1^a}+
\chi(\theta_1^a).
\end{equation}
\end{corollary}

\subsection{The Jacobian of the exponential map for a symmetric pair.}
\label{sec:Berezinians10} In this subsection, we assume that $\q$
is free of finite rank and we compute the Jacobian of the
exponential map $\q \to G/H$, as defined by formula
(\ref{eq:jacobian5}).

Let $(e_i)_{1\leq i \leq n}$ be a basis of $\q$, and
$(x^i)_{1\leq i \leq n}$ be
the dual basis of $\q^*$. Then $F(\q)=\K[[x^1,\dots,x^n]]$. The
generic point $y$ of $\q$ is equal to $y=\sum e_i x^i\in
F(\q,\q)=\q\otimes F(\q)$. Since $F(\q,\g)=F(\q,\h)\oplus
F(\q,\q)$, we consider $y$ as an element of the Lie
$F(\q)$-superalgebra $F(\q,\g)$.

If $k \in 2\N$, $\ad^k (y)$ stabilizes $\q\otimes F(\q)$ and we
denote by $\str_\q(\ad^k (y)) \in F(\q)_\0$ its supertrace, which
is an homogeneous polynomial of degree $k$ in the $x^i$. More
generally, for $w(t)=w_0+ w_2 t^2+\cdots \in \K_\0[[t^2]]$,
$\str_{\q}(w(\ad y))$ is defined and we have the Taylor expansion
\begin{equation}
\label{eq:adyk222}
  \str_{\q}(w(\ad  y))
= \sum_{k\geq 0}
  w_{k}
  \str_{\q}(\ad^{k} (y)) \in
F(\q)_\0.
\end{equation}
Similarly,   for $r(t)=1+ r_2 t^2+\cdots \in \K_\0[[t^2]]$, $r(\ad
y)$  induces an even invertible $F(\q)$-linear transformation of
$\q\otimes F(\q)$. We denote by $\Ber_\q(r(\ad y))$ its
Berezinian.

For $c\in \K_\0^*$ we define
\begin{equation}\label{eq:jacobian50}
J_c:=\Ber_{\q}\left(\frac{\sh(\ad  \frac{y}{c})}{\ad
\frac{y}{c}}\right) \in F(\q)_\0,
\end{equation}
and
\begin{equation}\label{eq:jacobian51}
w_c(t)=\log\left(\frac{\sh(\frac{t}{c})}{ \frac{t}{c}}\right) \in
t^2\K_\0[[t^2]].
\end{equation}
Because of formula (\ref{eq:differentielle4}) and
\begin{equation}\label{eq:jacobian52}
  w'_c(t)
= \frac{1}{c\, t}
  \left( t \coth\left( \frac{t}{c} \right)- c
  \right)
= \frac{1}{p_c(0) }\frac{p_c(t)-p_c(0)}{t},
\end{equation}
we have
\begin{equation}\label{eq:jacobian53}
w_c(t)=\sum_{n=1}^\infty b_{2n} \frac{
2^{2n}t^{2n}}{c^{2n}2n(2n)!}= \frac{1}{6 c^2}t^2 - \frac{1}{180
c^4}t^4 +\cdots
\end{equation}
Because of (\ref{eq:berk}) we obtain
\begin{equation}\label{eq:berradxc}
J_c =\exp\left(\str_\q(w_c(\ad y))\right).
\end{equation}
We write the beginning of the resulting Taylor expansion of $J_c$:
\begin{equation}\label{eq:berradxcc}
J_c =1+\frac{1}{6 c^2}\str_{\q}(\ad^{2} (y))
 - \frac{1}{180 c^4} \str_{\q}(\ad^4  (y))
 + \frac{1}{72 c^4}\str_{\q}^2(\ad^{2} (y))
 +\cdots
\end{equation}

\bigskip
We start with some computations on supertraces needed in the proof
of theorem \ref{thm:main}. For  $a\in \q$ and $i,j\in\N$ such that
$i+j\in 2\N+1$, the element $\ad^i (y) \ad a \ad^j( y)$ of
$\End_{F(\q)}(\g\otimes F(\q))$ stabilizes $\q\otimes F(\q)$.  Let
$t$ and $u$ be even formal variables. We use the notation
\begin{equation}\label{eq:Stu}
   S(t^iu^j:a)
= \str_\q
   \left(\ad^i (y)
\ad a  \ad^j( y)
   \right).
\end{equation}
Let  $k(t,u)=\sum_{i,j\geq 0} k_{ij}t^iu^j\in  \K[[t,u]]$ be such
that $k_{ij}=0$ if $i+j\in 2\N$. We define the locally finite sum
\begin{equation}\label{eq:Stuk}
  S(k(t,u):a)
= \sum_{i,j\geq 0}
  k_{ij}
  S(t^iu^j:a)
\in F(\q) .
\end{equation}
For example, for $q(t)\in t\K_\0[[t^2]]$ we have $q(\ad y)(a)\in
 F(\q,\h)$, and $\ad( q(\ad y)(a)) $ stabilizes $\q\otimes F(\q)$. We have
\begin{equation}\label{eq:adStuk}
\str_\q(\ad (q(\ad y)(a)))= S(q(t-u):a).
\end{equation}

\begin{lemma}\label{lem:Stu}
$S((t^2-u^2)k(t,u):a)=0$.
\end{lemma}
\begin{proof}
It is sufficient to prove it for a monomial $k=t^iu^j$. Let $Y$
and $X$ be the restrictions of $  \ad^2 (y) $ and $\ad^i( y) \ad a
\ad^j (y)$ to $\q\otimes F(\q)$. Since $Y$ is even, we have
$\str_\q(Y X)=\str_\q(X  Y)$ which proves the lemma.
\end{proof}

\medskip
Let $g(t)=g_0+g_2 t^2+ \cdots \in \K_0[[t^2]]$. We consider the
vector field $g(\ad y)(a)$ on $\q $.
\begin{proposition}\label{prop:divsym}
We have
$$
  \diverg(\zeta_{g(\ad y)(a)})
 = S \left(-\frac{g(t)-g(t-u)}{u}: a
     \right)
  .
$$
\end{proposition}
\begin{proof}
It is sufficient to  prove the statement for a monomial $g(t)=t^n$
with $n$ even. We argue like in the proof of proposition
\ref{prop:div}, replacing $x$ with $y$, and $\str$ with $\str_\q$.
Proposition \ref{prop:divsym} is equivalent to formula
(\ref{eq:divergenced}).
\end{proof}

\medskip
Let $r(t)=r_0+r_2 t^2+ \cdots \in \K_0[[t^2]]$. We consider the
function $\str_\q(r(\ad y))\in F(\q)$.
\begin{proposition}\label{prop:divsym2}
We have
$$
 \zeta_{g(\ad y)(a)}
 \left( \str_\q(r(\ad y))
 \right)=
 S \left( \frac{r(t)-r(u)}{t-u}g(t-u): a
   \right)
 .
$$
\end{proposition}
\begin{proof}
It is sufficient to  prove the statement for a monomial $r(t)=t^n$
with $n$ even. We argue like in the proof of proposition
\ref{prop:diffbera}, replacing $x$ with $y$, and $\str$ with
$\str_\q$. Proposition \ref{prop:divsym} is equivalent to formula
(\ref{eq:divergencedd}).
\end{proof}

\begin{lemma}\label{lem:calculcalcul}
Let $c \in \K_\0^\times$. For all $a\in \g$ we have
\begin{equation}\label{eq:toprove}
\zeta_{\alpha_c^a}\left({\str_\q(w_c (\ad y))}\right)+
\diverg{\zeta_{\alpha_c^a}}-\str_\q(\theta_c^a)=0.
\end{equation}
\end{lemma}
\begin{proof}
 We first consider the case $a\in \q$. From the definition
(\ref{eq:vectorfield}) of $\alpha_c^a$ and propositions
\ref{prop:divsym2} and \ref{prop:divsym}, we obtain
\begin{equation}\label{eq:jacobian445}
  \zeta_{\alpha_c^a}\left({\str_\q(w_c(\ad y))}\right)
= S \left( \frac{w_c (t)-w_c (u)}{t-u}p_c(t-u):a
    \right).
\end{equation}
and
\begin{equation}
\label{eq:jacobian446}
 \diverg{\zeta_{\alpha_c^a}}
= S \left(-\frac{p_c(t)-p_c(t-u)}{u}: a \right).
\end{equation}
From the definition (\ref{eq:thetac}) of $\theta_c^a$  and
equation (\ref{eq:adStuk}) we obtain
\begin{equation}\label{eq:jacobian447}
-\str_\q(\theta_c^a)= -S(q_c(t-u):a).
\end{equation}
Thus, the left hand side of (\ref{eq:toprove}) is equal to
$S(s(t,u):a)$ with
$$
s(t,u):=
\frac{w_c(t)-w_c(u)}{t-u}p_c(t-u)-\frac{p_c(t)-p_c(t-u)}{u}-q_c(t-u).
$$
Lemma \ref{lem:Stu} shows that to prove (\ref{eq:toprove}), it is
sufficient to prove that $s(t,u)$ belongs to the ideal
$(t^2-u^2)\K_\0[[t,u]]$.  For this, it is sufficient to prove that
we have $s(t,t)=0$ and $s(t,-t)=0$. We have (see
(\ref{eq:jacobian53}))
\begin{equation}\label{eq:jacobian448}
s(t,t)= w_c'(t)p_c(0) -\frac{p_c(t)-p_c(0)}{t} \mbox{ and }
s(t,-t)=\frac{p_c(t)-p_c(2t)}{t}-q_c(2t).
\end{equation}
Strangely enough, we already noticed the relation $s(t,t)=0$  in
formula (\ref{eq:jacobian52}) and the relation $s(t,-t)=0$ in
(\ref{eq:qccc}).

Last, we consider the simpler case  $a\in \h$. The proof is the
same with $p_c(t)$ replaced by $-t$, and $q_c(t)$ replaced by $1$.
We have $s(t,t)=-\frac{-t}{t}-1=0$ and $s(t,-t)=\frac{-t +
2t}{t}-1=0$.
\end{proof}

\begin{theorem} \label{thm:main}
Let $\g=\h\oplus \q$ be a symmetric pair such that $\q$ is free of
finite rank. The Jacobian $J\in F(\q)$ of the exponential map $\q
\to G/H$ is given by the formula $J=J_1$.
\end{theorem}
\begin{proof}
We verify that $J_1$ satisfies the conditions of lemma
\ref{lem:calcul}. We already showed that   $\delta(J_1)=1$ in
formula (\ref{eq:berradxcc}). Let $c\in \K_\0^*$. For a vector
field $\alpha \in F(\q,\q)$ we  get from (\ref{eq:berk}) the
relation
$$
  \frac{\zeta_\alpha(J_c)}
       {J_c}
= \zeta_\alpha
  \left({\str_\q(w_c(\ad y))}
  \right).
  $$   Let $a\in \g$. We recall  the
notations $\alpha_c^a\in F(\q,\q)$ and $\theta_c^a \in F(\q,\h)$
in theorem \ref{theo:formula}. The second condition
(\ref{eq:module3}) of lemma \ref{lem:calcul} follows   from lemma
\ref{lem:calculcalcul} choosing $c=1$.
\end{proof}

\subsection{Gorelik elements.}\label{sec:gorelik}
We recall  the definition (\ref{eq:twisted}) of the twisted
adjoint representation $\ad'$ of $\g$ in $U(\g)$. We extend $\ad'$
to a representation of $U(\g)$ in $U(\g)$ and we consider the map
$\tilde \gamma : U(\g) \to U(\g)$ defined by
\begin{equation}\label{eq:gamma1}
\tilde \gamma(u)=\ad'(u)(1) \mbox{ for } u\in U(\g).
\end{equation}
\begin{remark}
\label{rem:mia} By definition the map $\tilde\gamma$ intertwines
the left regular representation  and the twisted adjoint action.
\end{remark}
Moreover, since $\ad'(a)(1)=\ad(a)(1)=0$ for $a\in \h$, the map
$\tilde \gamma$ gives a map
\begin{equation}\label{eq:gamma2}
 \gamma : U(\g)/U(\g) \h \to U(\g).
\end{equation}
\begin{lemma}{\bf (Gorelik \cite{Gorelik2001}).}
\label{lem:gamma4}
  We have $ \gamma\circ \beta|_{S(\q)}=\beta\circ I_2$.
\end{lemma}
For example, for $b\in \q$, we have $\gamma(b)=\ad'(b)(1)=b  1-1
\sigma(b)=b+b=2b$. Thus  $\gamma$ is an isomorphism of coalgebras
from $U(\g)/U(\g) \h$ to $\beta(S(\q))$.
\begin{proposition}\label{prop:gamma}
\begin{enumerate}
    \item[i)]
The submodule $\beta(S(\q))\subset U(\g)$ is stable  under the
twisted adjoint action;
    \item[ii)]
as a $\g$-module it is isomorphic to $ U(\g)/U(\g) \h $;
    \item[iii)]
the map $\beta : S(\q) \to \beta(S(\q))$ intertwines the
representation $C_2$ with the twisted adjoint action.
\end{enumerate}
\end{proposition}
\begin{remark}
The first assertion is due to \cite{ABF1997}, the second to
\cite{Gorelik2001}, and the last one follows from remark
\ref{rem:mia}, lemma \ref{lem:gamma4}, theorem
\ref{theorem:emanuela3}, and remark \ref{rem:equivalent}.
\end{remark}

\bigskip

We assume moreover that the dimension of $\q$ is $(0,q)$ and that
the unimodularity condition(\ref{eq:unimb}) holds.  Let
$(e_i)_{1\leq i\leq q}$ be a totally ordered basis of $\q$, and
$(x^i)$ the dual basis of $\q^*$. We have $F(\q)=
\K[x^1,\dots,x^q]$ and the formal sum (\ref{eq:berradxcc}) which
defines $J_c$ is finite. We recall that, by theorems
\ref{thm:correspondance} and \ref{thm:main}, the element
$$
\widetilde \beta(J_1 e_1\dots e_q) \in U(\g)/U(\g)\h
$$
is basis of the $\K$-module of $\g$-invariant elements of
$U(\g)/U(\g)\h$, and a basis of $U(\g)/U(\g)\h$ as a
$F(G/H)$-module.

\begin{theorem}\label{thm:Gorelikfinal}
Let $\g=\h\oplus \q$ be a symmetric pair such that $\q$ is free of
finite rank $(0,q)$ and such that the unimodularity condition
(\ref{eq:unimb}) holds.   Let $(e_i)_{1\leq i\leq q}$ be a totally
ordered basis of $\q$. The element $\beta(J_2 e_1\dots e_q) \in
U(\g)$ is a basis of the $\K$-module of the invariants of the
twisted adjoint action of $\g$ in $\beta(S(\q))$.
\end{theorem}
\begin{proof}
Recall that $J_1 e_1\dots e_q $ is an invariant  for the
representation $\g\ni a\to C_1^a$ (see
theorem~\ref{theorem:emanuela3}, ii). Then the theorem follows
from proposition \ref{prop:gamma} part iii,  remark
\ref{rem:equivalent}, and $I_c(J_1  e_1\dots e_q)= c^q J_c
e_1\dots e_q$.
\end{proof}

The following corollary is the particular case of theorem
\ref{thm:Gorelikfinal} corresponding to $\K=\K_\0$, $\h=\g_\0$,
 and $\q= \g_\1$  free of finite rank. We state it in the
conditions of the introduction.

\begin{corollary}\label{cor:Gorelikfinal}
Let $\g=\g_\0\oplus \g_\1$ be a Lie superalgebra of finite
dimension $(p,q)$ over a field $\K$ of characteristic $0$, and
 $d \in S^q(\g_\1)$ with  $d\neq 0$.
Assume that the unimodularity assumption (\ref{eq:unim}) holds.
Then $\beta( J_2 d)$ is a Gorelik element of $U(\g)$.
\end{corollary}


\begin{bibdiv}
\begin{biblist}


\bib{ABF1997}{article}{
   author={Arnaudon, D.},
   author={Bauer, M.},
   author={Frappat, L.},
   title={On Casimir's ghost},
   journal={Comm. Math. Phys.},
   volume={187},
   date={1997},
   number={2},
   pages={429--439},
   issn={0010-3616},
}


\bib{BellFarnsteiner1993}{article}{
   author={Bell, Allen D.},
   author={Farnsteiner, Rolf},
   title={On the theory of Frobenius extensions and its application to Lie
   superalgebras},
   journal={Trans. Amer. Math. Soc.},
   volume={335},
   date={1993},
   number={1},
   pages={407--424},
   issn={0002-9947},
}



\bib{Berezin1970}{article}{
   author={Berezin, F. A.},
   author={Kac, G. I.},
   title={Lie groups with commuting and anticommuting parameters},
   language={Russian},
   journal={Mat. Sb. (N.S.)},
   volume={82 (124)},
   date={1970},
   pages={343--359},
}

\bib{Blattner1969}{article}{
   author={Blattner, Robert J.},
   title={Induced and produced representations of Lie algebras},
   journal={Trans. Amer. Math. Soc.},
   volume={144},
   date={1969},
   pages={457--474},
   issn={0002-9947},
}
\bib{Chemla1993}{article}{
   author={Chemla, Sophie},
   title={Cohomologie locale de Grothendieck et repr\'esentations induites
   de superalg\`ebres de Lie},
   language={French},
   journal={Math. Ann.},
   volume={297},
   date={1993},
   number={2},
   pages={371--382},
   issn={0025-5831},
}
\bib{Chemla1994}{article}{
   author={Chemla, Sophie},
   title={Propri\'et\'es de dualit\'e dans les repr\'esentations coinduites
   de superalg\`ebres de Lie},
   language={French, with English and French summaries},
   journal={Ann. Inst. Fourier (Grenoble)},
   volume={44},
   date={1994},
   number={4},
   pages={1067--1090},
   issn={0373-0956},
}
\bib{Chemla1994b}{article}{
   author={Chemla, Sophie},
   title={Poincar\'e duality for $k$-$A$ Lie superalgebras},
   journal={Bull. Soc. Math. France},
   volume={122},
   date={1994},
   number={3},
   pages={371--397},
   issn={0037-9484},
}


\bib{Cohn1963}{article}{
   author={Cohn, P. M.},
   title={A remark on the Birkhoff-Witt theorem},
   journal={J. London Math. Soc.},
   volume={38},
   date={1963},
   pages={197--203},
   issn={0024-6107},
}

\bib{Gorelik2001}{article}{
   author={Gorelik, Maria},
   title={On the ghost centre of Lie superalgebras},
   journal={Ann. Inst. Fourier (Grenoble)},
   volume={50},
   date={2000},
   number={6},
   pages={1745--1764 (2001)},
   issn={0373-0956},
}


\bib{Gorelik2002}{article}{
   author={Gorelik, Maria},
   title={Strongly typical representations of the basic classical Lie
   superalgebras},
   journal={J. Amer. Math. Soc.},
   volume={15},
   date={2002},
   number={1},
   pages={167--184 (electronic)},
   issn={0894-0347},
}


\bib{GL1999}{article}{
   author={Gorelik, Maria},
   author={Lanzmann, Emmanuel},
   title={The annihilation theorem for the completely reducible Lie
   superalgebras},
   journal={Invent. Math.},
   volume={137},
   date={1999},
   number={3},
   pages={651--680},
   issn={0020-9910},
}


\bib{Kosmann2000}{article}{
  author={Kosmann-Schwarzbach, Yvette},
  author={Monterde, Juan},
  title={Divergence operators and odd Poisson brackets},
  journal={Ann. Inst. Fourier},
  volume={52},
  date={2002},
  number={2},
  pages={419--456},
}


\bib{Koszul1982}{article}{
   author={Koszul, J.-L.},
   title={Graded manifolds and graded Lie algebras},
   conference={
      title={},
      address={Florence},
      date={1982},
   },
   book={
      publisher={Pitagora},
      place={Bologna},
   },
   date={1983},
   pages={71--84},
}



\bib{Lesniewski1995}{article}{
   author={Lesniewski, Andrzej},
   title={A remark on the Casimir elements of Lie superalgebras and
   quantized Lie superalgebras},
   journal={J. Math. Phys.},
   volume={36},
   date={1995},
   number={3},
   pages={1457--1461},
   issn={0022-2488},
}


\bib{Manin1997}{book}{
   author={Manin, Yuri I.},
   title={Gauge field theory and complex geometry},
   series={Grundlehren der Mathematischen Wissenschaften [Fundamental
   Principles of Mathematical Sciences]},
   volume={289},
   edition={2},
   publisher={Springer-Verlag},
   place={Berlin},
   date={1997},
   pages={xii+346},
   isbn={3-540-61378-1},
}

\bib{Musson1997}{article}{
   author={Musson, Ian M.},
   title={On the center of the enveloping algebra of a classical simple Lie
   superalgebra},
   journal={J. Algebra},
   volume={193},
   date={1997},
   number={1},
   pages={75--101},
   issn={0021-8693},
}


\bib{Petracci-thesis}{thesis}{
   author={Petracci, Emanuela},
   title={Functional equations and Lie algebras},
   language={English},
   organization={Tesi di Dottorato, Universit{\`a} di Roma ``La
   Sapienza''},
   date={2003},
}

\bib{Petracci2003}{article}{
   author={Petracci, Emanuela},
   title={Universal representations of Lie algebras by coderivations},
   journal={Bull. Sci. Math.},
   volume={127},
   date={2003},
   number={5},
   pages={439--465},
   issn={0007-4497},
}


\bib{Pinczon1990}{article}{
   author={Pinczon, Georges},
   title={The enveloping algebra of the Lie superalgebra ${\rm osp}(1,2)$},
   journal={J. Algebra},
   volume={132},
   date={1990},
   number={1},
   pages={219--242},
   issn={0021-8693},
}


\bib{Radford1986}{article}{
   author={Radford, David E.},
   title={Divided power structures on Hopf algebras and embedding Lie
   algebras into special-derivation algebras},
   journal={J. Algebra},
   volume={98},
   date={1986},
   number={1},
   pages={143--170},
   issn={0021-8693},
}



\bib{Rouviere1986}{article}{
   author={Rouvi{\`e}re, Fran{\c{c}}ois},
   title={Espaces sym\'etriques et m\'ethode de Kashiwara-Vergne},
   language={French, with English summary},
   journal={Ann. Sci. \'Ecole Norm. Sup. (4)},
   volume={19},
   date={1986},
   number={4},
   pages={553--581},
   issn={0012-9593},
}



\bib{Rouviere1994}{article}{
   author={Rouvi{\`e}re, Fran{\c{c}}ois},
   title={Fibr\'es en droites sur un espace sym\'etrique et analyse
   invariante},
   language={French, with English and French summaries},
   journal={J. Funct. Anal.},
   volume={124},
   date={1994},
   number={2},
   pages={263--291},
   issn={0022-1236},
}


\bib{Scheunert1979}{book}{
   author={Scheunert, Manfred},
   title={The theory of Lie superalgebras},
   series={Lecture Notes in Mathematics},
   volume={716},
   publisher={Springer},
   place={Berlin},
   date={1979},
   pages={x+271},
   isbn={3-540-09256-0},
}

\bib{Voronov1991}{book}{
   author={Voronov, T.},
   title={Geometric integration theory on supermanifolds},
   series={Soviet Scientific Reviews, Section C: Mathematical Physics Reviews},
   volume={9},
   publisher={Harwood Academic Publishers},
   place={Chur},
   date={1991},
   pages={iv+138},
   isbn={3-7186-5199-8},
}

\bib{MPV1990}{article}{
   author={Voronov, Alexander A.},
   author={Manin, Yu. I.},
   author={Penkov, I. B.},
   title={Elements of supergeometry},
   language={Russian},
   note={Translated in J. Soviet Math.\ {\bf 51} (1990), no.\ 1,
   2069--2083},
   conference={
      title={Current problems in mathematics. Newest results, Vol.\ 32},
   },
   book={
      series={Itogi Nauki i Tekhniki},
      publisher={Akad. Nauk SSSR Vsesoyuz. Inst. Nauchn. i Tekhn. Inform.},
      place={Moscow},
   },
   date={1988},
   pages={3--25},
}

\end{biblist}
\end{bibdiv}
\end{document}